    \documentclass[12pt]{article}
    \usepackage{amsmath}
    \usepackage{amssymb}
    \usepackage{amsthm}

\usepackage{pstricks}

\usepackage{psfig}

      \newcommand {\al}   {\alpha}          \newcommand {\bt}  {\beta}
                
      \newcommand {\del}  {\delta}          \newcommand {\Del} {\Delta}
              \newcommand {\ve}   {\varepsilon}
                 \newcommand {\vphi} {\varphi}
      \newcommand {\lam}  {\lambda}         \newcommand {\Lam}  {\Lambda}
      
                \newcommand {\Om}  {\Omega}
      \newcommand {\pl}   {\partial}        \newcommand {\s}    {\sigma}

      \newcommand {\RRR}  {{\mathbb R}}     \newcommand {\AAA}  {{\cal A}}
              \newcommand {\MMM}  {{\cal M}}
           \newcommand {\NNN}  {{\mathbb N}}
              \newcommand {\BBB}  {{\cal B}}
      \newcommand {\FFF}  {{\cal F}}          \newcommand {\SSS}  {{\cal S}}
      \newcommand {\QQQ}  {{\cal Q}}        
      \newcommand {\TTT}  {{\cal T}}
    \newcommand {\bbb}{B}      
      \newcommand {\interval}  {[-\pi/2,\, \pi/2]}

      \newtheorem{theor}{Theorem}
      \newtheorem{lem}{Lemma}
      
      \newtheorem{zam}{Note}
      \newtheorem{corollary}{Corollary}

\author{A. Yu. Plakhov\thanks{Department of Mathematics, Aveiro
University, Aveiro 3810, Portugal}}

\title{Billiards and two-dimensional problems\\ of optimal
resistance}

\date{}

\begin{document}

\maketitle

\begin{abstract}
A body moves in a medium composed of noninteracting point particles;
interaction of particles with the body is absolutely elastic. It is
required to find the body's shape minimizing or maximizing
resistance of the medium to its motion. This is the general setting
of optimal resistance problem going back to Newton.

Here, we restrict ourselves to the two-dimensional problems for
rotating (generally non-convex) bodies. The main results of the
paper are the following. First, to any compact connected set with
piecewise smooth boundary $B \subset \RRR^2$ we assign a measure
$\nu_B$ on $\pl(\text{conv}B) \times \interval$ generated by the
billiard in $\RRR^2 \setminus B$ and characterize the set of
measures $\{ \nu_B \}$. Second, using this characterization, we
solve various problems of minimal and maximal resistance of rotating
bodies by reducing them to special Monge-Kantorovich problems.
\end{abstract}

\begin{quote}
\end{quote}

\begin{quote}
\end{quote}

\section{Introduction}

A body moves in a homogeneous medium consisting of point particles.
The medium is very rare, so that mutual interaction of particles can
be neglected. The particles interact with the body in the absolutely
elastic manner; each particle performs several (maybe zero)
collisions and moves freely afterwards. It is required to find a
body, out of a given class of bodies, such that resistance of the
medium to the body's motion is minimal.

This is the general setting of minimal resistance problem. In order
to specify it, one needs to determine the class of bodies, the kind
of motion (for example, translational motion or a combination of
translational and rotational motion), and the state of medium (the
particles may rest as well as perform thermal motion). The problem
can be stated not only in $\RRR^3$, but also in spaces of other
dimensions; besides, the problems of {\it maximal} resistance can be
considered as well.

This problem goes back to Newton; he considered the class of convex
axially symmetric bodies of fixed length along the symmetry axis and
fixed maximal cross section orthogonal to this axis. A body of the
class moves in a medium of particles at rest, the velocity being
parallel to the symmetry axis. The solution of this problem is a
body whose boundary is composed of front and rear flat disks and a
smooth strictly convex lateral surface \cite{N}.

Since then, the problem has been studied by many mathematicians,
including Euler and Legendre; among other results, there were found
solutions in various classes of axisymmetric convex bodies. Since
the early 1990th, there have been obtained new interesting results
related to the minimization problem in classes of non-symmetric
and/or non-convex bodies \cite{BK}-\cite{LO}. It was proved, in
particular, that in the class of convex (generally non-symmetric)
bodies of fixed length and width the solution exists and does not
coincide with Newton's one \cite{BK,BrFK,BFK,LP2,BG}. The solution
was obtained numerically \cite{LO}. There was found analytically the
solution in the more restricted class of convex bodies with
developable lateral surface \cite{LP1}. Also, there were obtained
solutions in some classes of non-convex bodies satisfying the
so-called {\it single impact assumption} meaning that each particle
collides with the body at most once \cite{CL1,CL2}. Further, it was
shown that generically, in classes of non-convex bodies where
multiple collisions are allowed, infimum of resistance is zero, that
is, there exist almost perfectly streamlined bodies \cite{P1,P2}.
This result is based on the study of some special billiards in
unbounded regions in $\RRR^3$. It is essentially three-dimensional:
for the two-dimensional analogue of the problem, infimum is positive
\cite{P2}.

The problems involving rotational motion of bodies and/or thermal
motion of particles are more relevant to real life. Some of them,
concerning classes of {\it convex} bodies, are solved in
\cite{sb-math-averaged04} and \cite{temperature}, for the case of
slowly uniformly rotating convex bodies of fixed volume and the case
of translational motion of axially symmetric convex bodies in media
of positive temperature, respectively. As applied to classes of {\it
non-convex} bodies, these problems seem to be more difficult: to
calculate resistance, one needs to know the relation between the
initial and final velocity of each particle interacting with the
body; in the case of multiple collisions the calculation may turn
out to be very involved. Nevertheless, as will be shown below, these
difficulties can be overcome. Here, we restrict ourselves to the
two-dimensional case; the three-dimensional case is postponed to the
future.

In what follows, instead of a body moving in a medium, we shall
consider a flux of particles falling on the body. The body may rest
and may rotate around a fixed point. This picture is equivalent to
the initial one and usually is more convenient.

Our approach is as follows. Each body (a compact connected subset of
$\RRR^2$ with piecewise smooth boundary) is limited by a curve
composed of a ''convex part'' and a number of ''cavities''. Each
particle, interacting with the body, either reflects only once from
the convex part of curve, or, otherwise, gets into a cavity, makes
there a series of reflections, and eventually gets out of the cavity
and leaves the body forever. Define the angles of ''getting in'' and
''getting out'' by $\vphi$ and $\vphi^+$; to any cavity one assigns
a measure describing the joint distribution of $\vphi$ and
$\vphi^+$. Resistance of the body is determined by these {\it
measures generated by the cavities}, therefore characterization of
these measures is a key question for a wide range of problems of
minimal and maximal resistance.

This characterization is the main issue of the present work: we
determine closure, in the weak topology, of the set of measures
generated by cavities. The proof of this result is based on
construction of a family of cavities of special form and on a
detailed analysis of billiard dynamics in such cavities. We then
apply this result to the problem of optimal {\it mean resistance}
for the classes of (generally non-convex) bodies containing or being
contained in a given convex bounded set $K$. The bodies are subject
to slow uniform motion. By mean resistance, we understand the time
averaged value of resistance. The problem reduces to some special
one-dimensional Monge-Kantorovich problems of optimal mass transfer;
solving them, one finds that infimum and supremum of mean resistance
are equal to 0.9878... and to 1.5, respectively, where resistance of
$K$ is taken to be 1.

The paper is organized as follows. In section 2, to each body $B$
one assigns a measure $\nu_B$: the linear combination of a measure
generated by the convex part of its boundary and measures generated
by the cavities. The main theorem, consisting in characterization of
the set of measures $\nu_B$, is stated and then applied to the
problem of optimal mean resistance for rotating bodies. In section 3
the auxiliary result, consisting in characterization of the set of
measures generated by cavities, is formulated. Basing on this
result, we prove the main theorem. Finally, in section 4, the
auxiliary  result is proved.

\section{Main theorem}

\subsection{Statement of main theorem}

Let $\bbb \subset \RRR^2$ be a compact connected subset of Euclidean
space $\RRR^2$ with piecewise smooth boundary. The last means that
$\bbb$ is determined by a finite number of relations $f_i(x) \ge
0$,\, $f_i \in C^1(\RRR^2)$, besides $\nabla f_i(x) \ne 0$, whenever
$f_i(x) = 0$, and the vectors $\nabla f_i(x)$ and $\nabla f_j(x)$
are not collinear, whenever $f_i(x) = 0$,\, $f_j(x) = 0$,\, $i \ne
j$.

Consider the billiard in $\RRR^2 \setminus \bbb$. Consider a
billiard particle whose trajectory intersects the convex hull of
$\bbb$,\, conv$\bbb$. Initially the particle moves freely in $\RRR^2
\setminus \bbb$ with a unit velocity $v$. Let $x$ be the point of
first intersection of the particle with $\pl(\text{conv}\bbb)$;
denote by $n_x$ the unit outer normal vector to
$\pl(\text{conv}\bbb)$ at $x$ and denote by $\vphi \in \interval$
the angle between $n_x$ and $-v$. Let us agree that the angle is
counted from $n_x$ to $-v$ clockwise. Thus, to each particle motion
one assigns a value $(x,\vphi) \in \pl(\text{conv}\bbb) \times
\interval$.

After intersecting $\pl(\text{conv}\bbb)$, the particle moves inside
conv$\bbb \setminus \bbb$, elastically reflecting off the boundary
$\pl\bbb$, then intersects $\pl(\text{conv}\bbb)$ again and moves
freely afterwards. It may also happen that at some moment the
particle gets into a singular point of the boundary $\pl\bbb$, or
its trajectory touches the boundary, or it makes an infinite number
of reflections on a finite time interval, or the particle will stay
inside conv$\bbb$ forever. The set of values $(x,\vphi)$, for which
one of these events happens, has zero measure (see., e.g., \cite{T})
and therefore will be excluded from our consideration. Note that by
this agreement, the cases $\vphi = \pi/2$ and $-\pi/2$ are excluded,
therefore the particle intersects $\pl(\text{conv}\bbb)$ two times:
when getting {\it in} and when getting {\it out}.

Let $x^+ = x_\bbb^+(x,\vphi)$ be the point of second intersection.
Denote by $\vphi^+ = \vphi_\bbb^+(x,\vphi) \in \interval$ the angle
between the normal $n_{x^+}$ and the velocity of final free motion
$v^+$. Like $\vphi$, the angle $\vphi^+$ is measured from $n_{x^+}$
to $v^+$ clockwise.

Note that the set $\pl(\text{conv}\bbb) \setminus \pl\bbb$ is a
union of a finite or countable (maybe empty) family of mutually
disjoint open intervals, $\pl(\text{conv}\bbb) \setminus \pl\bbb =
I_1 \cup I_2 \cup \ldots$. The set conv$\bbb \setminus \bbb$ is
decomposed into several connected components; denote by $\Om_i$ the
closure of the connected component containing $I_i$,\, $i \ge 1$.
All the sets $\Om_i$ are different.
%\newpage
 
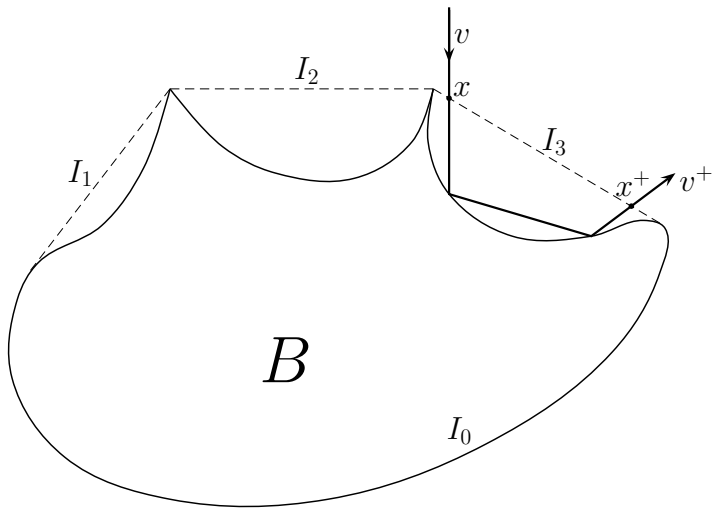
\vspace*{70mm}

\scalebox{0.7}{
 \pscurve(3,8)(4,6.9)(5,6.4)(6.5,6.3)(7.6,7)(8,8)
 \pscurve(3,8)(2.5,6.5)(1.7,5.4)(0.5,4.7)(0.1,4)(0,2.5)(1,1)(3,0.2)(6,0.2)(9,1.25)(11.4,3)(12.3,4.5)(12.4,5.35)(12,5.5)(11,5.2)(9.5,5.2)(8.3,6)(7.9,7)(8,8)
 \psline[linestyle=dashed,linewidth=0.4pt](3,8)(8,8)
 \psline[linestyle=dashed,linewidth=0.4pt](3,8)(0.35,4.55)
 \psline[linestyle=dashed,linewidth=0.4pt](12.3,5.45)(8,8)
 \psline[linewidth=1.2pt,arrows=->,arrowscale=1.5](8.3,9.5)(8.3,6)(11,5.2)(12.6,6.4)
 \psline[linewidth=1.2pt,arrows=->,arrowscale=1.5](8.3,9.55)(8.3,8.5)
 \rput(8.55,9){\Large $v$}
 \rput(13,6.3){\Large $v^+$}
 \rput(8.55,8){\Large $x$}
 \rput(11.8,6.2){\Large $x^+$}
 \psdots[dotsize=3pt](8.3,7.82)(11.76,5.77)
    \rput(8.5,1.5){\Large $I_0$}
    \rput(10.3,7){\Large $I_3$}
    \rput(5.6,8.35){\Large $I_2$}
    \rput(1.3,6.4){\Large $I_1$}
   }
   \rput(3.5,2){\Huge $\bbb$}

 \vspace{10mm}

Designate by $I_0 := \pl(\text{conv}\bbb) \cap \pl\bbb$ the {\it
convex part} of the boundary $\pl\bbb$; thus, one has
$$
\pl(\text{conv}\bbb) = I_0 \cup I_1 \cup I_2 \cup \ldots.
$$
If the point of first intersection $x$ belongs to $I_0$ then $x^+ =
x$,\, $\vphi^+ = -\vphi$, and the time interval of staying of the
particle within conv$\bbb$ reduces to a point. If, otherwise, $x$
belongs to an interval $I_i$,\, $i \ge 1$ then the intersection of
the particle trajectory with conv$\bbb$ entirely belongs to $\Om_i$
and $x^+$ belongs to the same interval $I_i$. In both cases one has
$n_{x^+} = n_x$.

Define a Borel measure $\mu = \mu_{\text{conv}\bbb}$ on the set
$\pl(\text{conv}\bbb) \times \interval$ by $d\mu(x,\vphi) =
\cos\vphi\, dx d\vphi$, where $dx$ is the element of length on
$\pl(\text{conv}\bbb)$. Consider the mapping $\TTT_\bbb:\, (x,\vphi)
\mapsto (x^+,\vphi^+)$. It is a one-to-one mapping of a full measure
subset of $\pl(\text{conv}\bbb) \times \interval$ onto itself;
moreover, the following holds:

{\bf T1}\, $\TTT_\bbb$ preserves the measure $\mu$;

{\bf T2} $\TTT_\bbb = \TTT_\bbb^{-1}$.\\
Actually, the conditions T1 and T2 are derived from the fact that a
billiard dynamical system preserves the Liouville measure and is
invariant with respect to time inversion.

The sets $I_i \times \interval$ are invariant with respect to the
mapping $\TTT_\bbb$; denote by $\TTT_\bbb^i$ the restriction of
$\TTT_\bbb$ to $I_i \times \interval$. The restriction of $\mu$ to
$I_i \times \interval$ will be also designated by $\mu$. Each
mapping $\TTT_\bbb^i$ transforms a subset of full measure of $I_i
\times \interval$ onto itself and also satisfies T1 and T2.

Denote by $|\pl(\text{conv}\bbb)|$ the length of the curve
$\pl(\text{conv}\bbb)$, by $|I_i|$, the length of $I_i$,\, $i = 1,\
2, \ldots$, and by $|I_0|$, the length of $I_0$:\, $|I_0| =
|\pl(\text{conv}\bbb)| - \sum_{i\ge 1} |I_i|$. Denote also $\kappa_i
= |I_i|/|\pl(\text{conv}\bbb)|$; one has $\sum_i \kappa_i = 1$.
Define Borel measures $\nu_\bbb$,\, $\nu_\bbb^i$ on the square $Q :=
\interval \times \interval$ as follows: for any measurable set $A
\subset Q$,
$$
\nu_\bbb(A) := \frac{1}{|\pl(\text{conv}\bbb)|}\, \mu \left( \left\{
(x,\vphi) \in \pl(\text{conv}\bbb) \times \interval:\,
(\vphi,\vphi_\bbb^+(x,\vphi)) \in A \right\} \right),
$$
$$
\nu_\bbb^i(A) := \frac{1}{|I_i|}\, \mu \left( \left\{ (x,\vphi) \in
I_i \times \interval:\, (\vphi,\vphi_\bbb^+(x,\vphi)) \in A \right\}
\right).
$$
These definitions imply that $\nu_\bbb = \sum_i \kappa_i
\nu_\bbb^i$.

Note that the measures $\nu_\bbb$,\, $\nu_\bbb^i$ remain unchanged
under translations, rotations and similarity transformations of
$\bbb$. The total measures are $\nu_\bbb(Q) = \nu_\bbb^i(Q) = 2$ for
any $i$.

Define the Borel measure $\lam$ on $\interval$ by $d\lam(\vphi) =
\cos\vphi\, d\vphi$, denote by $\pi^1$,\, $\pi^2:\, Q \to \interval$
the projections $\pi^1(\vphi, \vphi^+) = \vphi$,\, $\pi^2(\vphi,
\vphi^+) = \vphi^+$, and denote by $\pi^{\text diag}:\, Q \to Q$ the
symmetry with respect to the diagonal $\vphi = \vphi^+$:\,
$\pi^{\text diag}(\vphi, \vphi^+) = (\vphi^+, \vphi)$. Denote by
$\MMM$ the set of Borel measures $\nu$ on $Q$ such that

{\bf A1}\, both orthogonal projections of $\nu$ to the coordinate
axes coincide with $\lam$; that is, $\pi_\#^1 \nu = \lam = \pi_\#^2
\nu$;

{\bf A2}\, the measure $\nu$ is symmetric with respect to the
diagonal $\vphi = \vphi^+$; that is, $\pi_\#^{\text diag} \nu =
\nu$.

The measures $\nu_\bbb,\ \nu_B^i$ satisfy the conditions A1 and A2;
that is,
$$
\nu_\bbb,\ \nu_B^i \in \MMM.
$$
This property can be easily derived from the conditions T1 and T2
for the mappings $T_B$ and $T_B^i$; see also
\cite{sb-math-averaged04}.

Note that there exists a unique measure $\nu^0 \in \MMM$ whose
support belongs to the diagonal $\vphi = -\vphi^+$. If $|I_0|
> 0$ then one has $\nu_\bbb^0 = \nu^0$. Moreover, if $\bbb$
is convex then the family $\{ I_i \}$ contains a unique element
$I_0$, therefore $\nu_\bbb = \nu^0$.

Let us give a mechanical interpretation of the measures $\nu_\bbb$
and $\nu_\bbb^i$. A two-dimensional body $\bbb$ at rest is situated
in a medium of non-interacting identical point particles, moving in
all possible directions. The medium is uniform and isotropic, that
is, its density is constant and velocities of particles, contained
in any space region at any time instant, are uniformly distributed
on the unit circumference. When colliding with the body, the
particles are elastically reflected by it. For any particle,
interacting with the body, mark $\vphi$ and $\vphi^+$, the angles of
''getting in'' and ''getting out'' at the moments of intersection
with $\pl(\text{conv}\bbb)$. The measure $\nu_\bbb$ describes the
normalized joint distribution over $\vphi$ and $\vphi^+$ of the
number of particles that have been interacted with the body during a
fixed time interval. The normalizing factor is chosen in such a way
that the normalized total number of particles that have been
interacted with the body equals 2. The measures $\nu_\bbb^i$
describe the normalized joint distribution over $\vphi$ and
$\vphi^+$ of the number of particles that first intersected $I_i$
during a fixed time interval. The normalizing factor is chosen in
such a way that the normalized total number of particles that first
intersected $I_i$ during this time interval equals 2.

Any measure $\nu_\bbb$ belongs to $\MMM$. In the present paper we
are interested in an inverse question: {\it is it true that any
measure from $\MMM$ can be approximated by measures $\nu_\bbb$}? The
positive answer to this question is given by the following theorem.

Let $K_1$ and $K_2$ be compact convex sets such that $K_1 \subset
K_2$ and dist$(\pl K_1, \pl K_2) > 0$. Denote by $\BBB_{K_1,K_2}$
the class of compact connected sets $\bbb$ with piecewise smooth
boundary such that $K_1 \subset \bbb \subset K_2$.

\begin{theor}\label{t1}
The set $\{ \nu_\bbb, \ \bbb \in \BBB_{K_1,K_2} \}$ is everywhere
dense in $\MMM$ in the weak topology, that is, for any $\nu \in
\MMM$ there exists a family of sets $\{ \bbb_\ve,~ \ve
> 0 \} \subset \BBB_{K_1,K_2}$ such that for any continuous function $f :\, Q \to
\RRR$,
\begin{equation*}\label{weak convergence}
\lim_{\ve \to +0} \int_Q f\, d\nu_{\bbb_\ve} = \int_Q f\, d\nu.
\end{equation*}
\end{theor}

\subsection{Applications of main theorem}

Let us illustrate on examples that this theorem is useful when
solving problems of minimization and maximization of aerodynamic
resistance.
 \vspace{2mm}

{\bf Example 1}~ On a rotating body in $\RRR^2$ is incident a flux
of point particles. Density of the flux is constant; all the
particles have an equal velocity. When colliding with the body, the
particles are reflected by it according to the law of elastic
reflection. The particles do not mutually interact. The body is
fixed at some point and rotates around it with a constant angular
velocity. The rotation velocity is small: at each moment, the linear
velocity of any point of the body is much less than the flux
velocity. Thus, when considering interaction of the body with any
individual particle, one can neglect the effect of rotation.

The force of pressure of the flux on the body is a periodic
vector-valued function with the period equal to the period of the
body rotation. The mean value of this function --- let it call {\it
mean resistance} --- is a vector parallel to the flux velocity. It
is required to find infimum of mean resistance

(a) in the class of {\it convex} bodies of given area;

(b) in the class of (generally {\it non-convex}) bodies of given
area.

Denote by $\bbb$ the set occupied by the body at the zero instant.
The mean resistance is equal (up to a factor proportional to the
medium density and to the squared flux velocity) to
\begin{equation}\label{example1 mean R}
\bar R(\bbb) = |\pl (\text{conv}\bbb)| \cdot \FFF(\nu_\bbb),
\end{equation}
where
\begin{equation}\label{example1 mean resist}
\FFF(\nu) = \int\!\!\!\int_Q \left( 1 + \cos(\vphi - \vphi^+)
\right) d\nu(\vphi,\vphi^+).
\end{equation}

(a)~ If $\bbb$ is convex then $\nu_\bbb = \nu^0$ and conv$\bbb =
\bbb$, hence $\bar R(\bbb) = |\pl\bbb| \cdot \FFF(\nu^0)$. Thus, the
minimal resistance problem reduces to an isoperimetric problem: find
a convex set $\bbb$ with fixed area $S$ and minimal perimeter
$\pl\bbb$. The solution is a circle $\bbb^{(r)}$ of radius $r =
\sqrt{S/\pi}$; minimal resistance equals $\bar R(\bbb^{(r)}) =
2\sqrt{\pi S} \cdot \FFF(\nu^0)$. Note that
$$
\FFF(\nu^0) = \int\!\!\!\int_Q \left( 1 + \cos(\vphi - \vphi^+)
\right) d\nu^0(\vphi,\vphi^+) = \int_{-\pi/2}^{\pi/2} \left( 1 +
\cos 2\vphi \right) \cos\vphi\, d\vphi = 8/3.
$$

 \vspace{1mm}

(b)~ In the non-convex case the problem is solved in three steps.

(i) Find
$$
\inf_{\nu\in\MMM} \FFF(\nu);
$$
this is a one-dimensional {\it Monge-Kantorovich mass transport
problem} with the cost function $c(\vphi,\vphi^+) = 1 + \cos(\vphi -
\vphi^+)$ and with both marginal measures equal to $\lam$. It was
solved in \cite{sb-math-averaged04}; the minimizing measure $\nu_*$
satisfies the relation
$$
\frac{\FFF(\nu_*)}{\FFF(\nu^0)}\, =\, 0.9878...\,.
$$

(ii) Using theorem \ref{t1}, by diagonal method choose a sequence of
sets $\bbb_n'$ such that $\bbb^{(r)} \subset \bbb_n' \subset
\bbb^{(r + \frac 1n)}$ and the sequence of measures $\nu_{\bbb_n'}$
weakly converges to $\nu_*$.

(iii) Make similarity transformations $\bbb_n = k_n \bbb_n'$ in such
a way that the areas of all obtained sets $\bbb_n$ are equal to $S$.
The obtained sequence of sets is a solution of the minimization
problem.

Indeed, whatever the set $B$ of area $S$, the length $\pl
(\text{conv}\bbb)$ is not less than the perimeter of a circle of
area $S$, that is, $|\pl (\text{conv}\bbb)| \ge |\pl \bbb^{(r)}| =
2\sqrt{\pi S}$; moreover, $\FFF(\nu_\bbb) \ge \FFF(\nu_*)$. Hence,
$\bar R(\bbb) \ge 2\sqrt{\pi S} \cdot \FFF(\nu_*)$.

On the other hand, one has $\lim_{n \to \infty} k_n = 1$ and
$|\pl\bbb^{(r)}| \le |\pl(\text{conv} \bbb_n')| \le |\pl\bbb^{(r +
1/n)}|$, hence $\lim_{n\to\infty} |\pl(\text{conv} \bbb_n)| =
|\pl\bbb^{(r)}|$. Further, one has $\lim_{n\to\infty}
\FFF(\nu_{\bbb_n}) = \FFF(\nu_*)$. This implies that
$\lim_{n\to\infty} \bar R(\bbb_n) = 2\sqrt{\pi S} \cdot
\FFF(\nu_*)$.
 \vspace{1mm}

The ratio of least resistance for {\it non-convex} bodies of fixed
area to least resistance for {\it convex} bodies of the same area
equals
$$
\frac{\lim_{n\to\infty} \bar R(\bbb_n)}{\bar R(\bbb^{(r)})}\, =\,
\frac{\FFF(\nu_*)}{\FFF(\nu^0)}\, =\, 0.9878...\,.
$$
Thus, the gain is approximately $1.22\%$.
 \vspace{2mm}

{\bf Example 2}~ Like in the previous example, consider the
functional of mean resistance $\bar R(\bbb)$ (\ref{example1 mean
R},\ref{example1 mean resist}). The following problems are under
consideration:

(a)~ find $\inf \bar R(\bbb)$ in the class of sets $\bbb$ containing
a given convex bounded set $K$;

(b)~ find $\sup \bar R(\bbb)$ in the class of sets $\bbb$ contained
in a given convex bounded set $K$ with nonempty interior.

The problem (a) is solved in the same way as the problem (b) of the
previous example; one constructs a sequence of sets $\bbb_n^*$ such
that the sequence $\nu_{\bbb_n^*}$ weakly converges to $\nu_*$ and,
moreover, $K \subset \bbb_n^* \subset (1 + 1/n)\, K$. Here and in
what follows, $\kappa\, K := \{ \kappa x:\, x \in K \}$ designates
the set obtained from $K$ by the homothety with center at the origin
and ratio $\kappa$. The sequence $\bbb_n^*$ minimizes resistance;
the following relations hold: $\lim_{n\to\infty}
|\pl(\text{conv}\bbb_n^*)| = |\pl K|$ and $\lim_{n\to\infty}
\FFF(\nu_{\bbb_n^*}) = \FFF(\nu_*)$. On the other hand, one has
$\bar R(K) = |\pl K| \cdot \FFF(\nu^0)$, hence
$$
\frac{\lim_{n\to\infty} \bar R(\bbb_n^*)}{\bar R(K)} = 0.9878...\,.
$$

Now, consider the problem (b). Let us first note that there exists a
unique measure $\nu^\star \in \MMM$ whose support belongs to the
diagonal $\vphi = \vphi^+$. One has $\FFF(\nu^\star) =
\int_{-\pi/2}^{\pi/2} 2\cos\vphi\, d\vphi = 4$. This measure
maximizes the functional $\FFF$:
$$
\max_{\nu\in\MMM} \FFF(\nu) = \FFF(\nu^\star) = 4.
$$
Further, one constructs a sequence of sets $\bbb_n^\star$ such that
$\nu_{\bbb_n^\star}$ weakly converges to $\nu^\star$ and $\left(1 -
1/n\right) K \subset \bbb_n^\star \subset K$. This sequence
maximizes resistance: one has $\lim_{n\to\infty}
|\pl(\text{conv}\bbb_n^\star)| = |\pl K|$,\, $\lim_{n\to\infty}
\FFF(\nu_{\bbb_n^\star}) = \FFF(\nu^\star)$, therefore
$$
\frac{\lim_{n\to\infty} \bar R(\bbb_n^\star)}{\bar R(K)} = 1.5.
$$

Thus, there are constructed two sequences of sets approximating $K$
and providing solutions for the problems of minimal and maximal
resistance. The limit value of resistance for the first sequence is
$1.22\%$ less than resistance of $K$, and for the second sequence,
$50\%$ more. Note in passing that boundaries of the sets of both
sequences approximate boundary of the limit set, $\pl K$, in $C^0$,
but not in $C^1$.
 \vspace{2mm}

Like in these examples, one can put other problems of resistance
minimization and maximization for two-dimensional bodies performing
translational and/or rotational motion in rarefied medium, and
reduce these problems to special one-dimensional problems of
Monge-Kantorovich mass transfer. In this approach, of importance are
the assumptions that the medium is rarefied enough, so that mutual
interaction of particles can be neglected, and that the collisions
of particles with the body are absolutely elastic. The particles may
rest, and may perform chaotic thermal motion. The body's rotation
may be uniform or non-uniform. The velocity of rotational motion may
be small as compared to the velocity of translational motion, and
may be comparable to it. Many problems of this kind reduce to
minimization or maximization problems for functionals similar to
(\ref{example1 mean resist}). The functionals are defined on $\MMM$,
and the cost function is determined by a particular problem under
consideration. The work in this direction is in progress. Another
important task is generalization of this approach to the
three-dimensional case.

\section{Auxiliary result}

In this section, an auxiliary statement, theorem \ref{t2}, is
formulated, and basing on this statement, theorem \ref{t1} is
derived.

\subsection{Statement of the auxiliary result}\label{auxiliary}

Let $\Om \subset \RRR^2$ be a compact set with piecewise smooth
boundary. Denote by $dx$ the element of length on $\pl\Om$ and
define the Borel measure $\mu = \mu_\Om$ on $\pl\Om \times
\interval$ by $d\mu(x,\vphi) = \cos\vphi\, dx d\vphi$. The billiard
in $\Om$ generates the mapping $T = T_\Om$ in the following way. Let
$x \in \pl\Om$ be a regular point of the boundary and let $-\pi/2 <
\vphi < \pi/2$. Find the unit vector $v$ such that the angle between
$-n_x$ and $v$ is $\vphi$; recall that $n_x$ designates the unit
outer normal vector to $\pl\Om$ at $x$ and that the angle is counted
from $-n_x$ to $v$ clockwise. From $x$, launch a ray in the
direction $v$, denote by $x'$ the point of first intersection of the
ray with $\pl\Om$ and, if $x'$ is a regular point of the boundary,
reflect the vector $v$ off the boundary according to the law of
absolutely elastic reflection. The reflected vector equals $v' = v -
2(v,\, n_{x'})\, n_{x'}$. Denote by $\vphi'$ the angle that $v'$
makes with $-n_{x'}$; the angle is counted from $-n_{x'}$ to $v'$
clockwise. By definition, $T(x,\vphi) = (x',\vphi')$. It is well
known that $T$ is a one-to-one correspondence between two full
measure subsets of $\pl\Om \times \interval$ and preserves the
measure $\mu$.

For further convenience, introduce the notation $T^n(x,\vphi) =:
(x_n(x,\vphi),\, \vphi_n(x,\vphi))$,\, $n = 1,\ 2,\ldots$. Let $I
\subset \pl\Om$ be a Borel subset of $\pl\Om$ of positive
one-dimensional Hausdorff measure, $|I| > 0$. For any $(x, \vphi)
\in I \times \interval$, denote by $n(x,\vphi) = n_{\Om,I}(x,\vphi)$
the smallest positive integer $n$ such that $x_n(x,\vphi) \in I$.
Define the mapping $\TTT_{\Om,I}$ by $\TTT_{\Om,I}(x,\vphi) :=
(x_{n(x,\vphi)} (x,\vphi),\, -\vphi_{n(x,\vphi)} (x,\vphi))$. In
other words, we launch a billiard particle from the point $x \in I$
at the angle $\vphi$ and wait until it reflects off $I$ again. Fix
the point of the second reflection $x^+$ and the angle $\vphi^+$ of
the particle's motion just before the reflection. By definition,
$\TTT_{\Om,I}(x,\vphi) = (x^+, \vphi^+)$.

The mapping $\TTT_{\Om,I}$ preserves the measure $\mu$ and, by
Poincar\'e recurrence theorem, represents a one-to-one
correspondence between two full measure subsets of $I \times
\interval$.

Denote $\TTT_{\Om,I}(x,\vphi) =: (x_{\Om,I}^+(x,\vphi),\,
\vphi_{\Om,I}^+(x,\vphi))$ (in what follows we shall, as a rule,
omit the subscripts of the functions $x^+,\ \vphi^+$). Define the
Borel measure $\nu_{\Om,I}$ on $Q = \interval \times \interval$ as
follows. For any measurable set $A \subset Q$, put
$$
\nu_{\Om,I}(A) = \frac{1}{|I|}\, \mu \left( \{ (x,\vphi):\, (\vphi,
\vphi_{\Om,I}^+(x,\vphi)) \in A \} \right).
$$

Recalling the definition of sets $\Om_i$,\, $I_i$,\, $i \ge 1$,
mappings $\TTT_\bbb^i$ and measures $\nu_\bbb^i$, one sees that
$\TTT_\bbb^i = \TTT_{\Om_i,I_i}$ and $\nu_\bbb^i = \nu_{\Om_i,I_i}$.

It is easy to see that the measure $\nu_{\Om,I}$ satisfies the
conditions A1 and A2 and hence belongs to $\MMM$. We are interested
in an inverse problem: show that any measure $\nu \in \MMM$ can be
approximated by measures of the form $\nu_{\Om,I}$.

Denote by $\SSS$ the set of pairs $(\Om, I)$ such that $\Om \subset
\RRR^2$ is a compact set with piecewise smooth boundary, $I \subset
\pl\Om$ is a line segment, and $\Om$ lies on one side of the
straight line containing $I$ (see the picture below). The auxiliary
result is as follows.

\begin{theor}\label{t2}
The set $\{ \nu_{\Om,I},\, (\Om,I) \in \SSS \}$ is everywhere dense
in $\MMM$ in the weak topology, that is, for any measure $\nu \in
\MMM$ there exists a family of pairs $\{ (\Om(r), I(r)),\, r \ge 1
\} \subset \SSS$ such that for every continuous function $f:\, Q \to
\RRR$,
$$
\lim_{r\to+\infty} \int_Q f\, d\nu_{\Om(r), I(r)} = \int_Q f\, d\nu.
$$
\end{theor}

Note that in \cite{P-unbounded regions} there was obtained a
particular case of this theorem; namely, it was proved that
$\nu^\star$ can be approximated by measures $\nu_{\Om,I},\, (\Om,I)
\in \SSS$. Recall that $\nu^\star \in \MMM$ and the support of
$\nu^\star$ belongs to the diagonal $\vphi = \vphi^+$. This result
allows one to solve the mean resistance {\it maximization} problem
for slowly rotating bodies. Moreover, in \cite{P-unbounded regions}
there was obtained a three-dimensional analogue of this result,
which allows one to solve the three-dimensional maximization
problem.

\subsection{Derivation of main theorem}

Here theorem 1 is derived from theorem 2.

Denote by $\SSS'$ the set of pairs $(\Om, I) \in \SSS$ such that for
some positive $a$ and $b$,
$$
I = [-a,\, a] \times \{ - b \} \ \ \ \text{and} \ \ \ \Om \cap
\left( \RRR \times (-\infty,\, 0) \right) = [-a,\, a] \times [-b,\,
0)
$$
(see the figure below). The following lemma holds.
 \vspace*{43mm}

       \scalebox{1.5}{
       \rput(-1.5,-0.9){
 \psline[linewidth=0.5pt](1,1)(4,1)
 \psline[linewidth=1.3pt](2,1)(4,1)
        \psline[linewidth=0.6pt,arrows=->,arrowscale=1.2](2.3,1)(2.3,1.6)
        \psline[linewidth=0.6pt](2.3,1.6)(2.3,2.7)
           \psline[linewidth=0.6pt,arrows=->,arrowscale=1.2](2.3,2.7)(1.95,1.51)
           \psline[linewidth=0.6pt](1.95,1.51)(1.8,1)
               \psline[linewidth=0.6pt,arrows=->,arrowscale=1.2](1.8,1)(1.5,2.0)
               \psline[linewidth=0.6pt](1.8,1)(1.35,2.5)
        \psline[linewidth=0.6pt,arrows=->,arrowscale=1.2](1.35,2.5)(3.23,1.3)
        \psline[linewidth=0.6pt](3.23,1.3)(3.7,1)
 \pscurve[linewidth=0.5pt](4,1)(4.7,1.2)(5.2,2)(5.1,3.0)(4.4,3.1)(3.3,2.5)(2.3,2.7)(1.5,2.6)(0.8,1.8)(1,1)
 \psdots[dotsize=2.5pt](2,1)(4,1)
   \rput(4.3,2.2){\Large $\Om$}
   \rput(3.0,0.5){\large $I$}
        \rput(3.7,0.85){\tiny $x^+$}
        \rput(3.8,1.2){\tiny $\vphi^+$}
        \rput(2.4,0.8){\tiny $x,\! \vphi$}
       }}
   \rput(2.0,-1.6){An example of a pair $(\Om, I) \in \SSS$}
             \rput(6.5,0){
  \psline(1,1)(2,1)(2,0)(4,0)(4,1)
      \psline[linewidth=2pt](2,0)(4,0)
  \pscurve(4,1)(4.7,1.2)(5.2,2)(5.1,3.0)(4.4,3.1)(3.3,2.5)(2.3,2.7)(1.5,2.6)(0.8,1.8)(1,1)
   \rput(2.7,1.7){\huge $\Om$}
   \rput(3.2,0.3){\large $I$}
    \rput(6.1,3.0){\huge $\tilde\Om$}
  \pspolygon[linestyle=dashed](1,1)(2,1)(2,0)(4,0)(4,1)(5.5,1)(5.5,3.5)(0.5,3.5)(0.5,1)
  \psline(2,1)(0.5,1)
  \psline(4,1)(5.5,1)
         \rput(3,-0.8){An example of a pair $(\Om, I) \in \SSS'$}
   }

 \vspace{25mm}

\begin{lem}\label{lem vspomogat}
For any $(\Om, I) \in \SSS$ there exists a sequence $\{ (\Om^{[n]},
I^{[n]}),\ n = 1,\, 2,\ldots \} \subset \SSS'$ such that the
sequence of measures $\nu_{\Om^{[n]},I^{[n]}}$ weakly converges to
$\nu_{\Om,I}$.
\end{lem}

The proof of this lemma is not difficult and is placed in Appendix.

It follows from theorem \ref{t2} and lemma \ref{lem vspomogat} that
the set of measures $\{ \nu_{\Om,I}, \ (\Om,I) \in \SSS' \}$ is
everywhere dense in $\MMM$ in the weak topology.

Let $K_1$ and $K_2$ be compact convex sets, $K_1 \subset K_2$,\,
dist$(\pl K_1, \pl K_2) > 0$, and let $(\Om, I) \in \SSS'$ be an
arbitrary pair. To prove theorem \ref{t1}, it suffices to show that
the measure $\nu_{\Om,I}$ is a weak limit of measures of the form
$\nu_\bbb$,\, $\bbb \in \BBB_{K_1,K_2}$.

Denote $\tilde\Om_{a,b,c} = ([-a,\, a] \times [-b,\, 0]) \cup
((-c,\, c) \times [0,\, c))$. Obviously, there exist positive values
$a,\, b,\, c$ such that $a \le c$, $I = [-a,\, a] \times \{ -b \}$,
and $\Om \subset \tilde\Om_{a,b,c}$. We shall use the shorthand
notation $\tilde\Om_{a,b,c} = \tilde\Om$. Later on, we shall also
use the sets $\tilde{\Om}^C = \tilde{\Om}^C_{a,b,c} := (-c,\, c)
\times [-b,\, c)$,\, $\tilde\Om^l = \tilde\Om^l_{a,b,c} := (-c,\,
-a) \times [-b,\, 0)$, and $\tilde\Om^r = \tilde\Om^r_{a,b,c} :=
(a,\, c) \times [-b,\, 0)$ (see the figure below). One has
$\tilde\Om^C = \tilde\Om \cup \tilde\Om^l \cup \tilde\Om^r$; the
sets $\tilde\Om$,\, $\tilde\Om^l$,\, $\tilde\Om^r$ are mutually
disjoint. We shall also designate by $V$ an isometry of general
form. Thus, the composition of a homothety with ratio $k$ and center
at the origin and an isometry $V$, applied to a set $B \in \RRR^2$,
has the form $VkB$. For brevity, sets of this form will be referred
to as {\it copies} of $B$.

Consider a convex polygon $K_0$ such that $K_1 \subset K_0 \subset
K_2$ and dist$(\pl K_0, \pl K_1) > 0$, and fix $\ve > 0$. Our goal
is to construct a set of mutually disjoint copies of $\Om$ contained
in $K_0 \setminus K_1$ in such a way that the corresponding copies
of $I$ are contained in the boundary $\pl K_0$ and fill it all,
except possibly for a set of common length at most $\ve$. Then the
set $B_\ve$, obtained from $K_0$ by set-theoretic subtracting that
set of copies of $\Om$, will be defined. Making $\ve$ arbitrarily
small, one will be able to make the measure $\nu_{B_\ve}$
arbitrarily close (in variation) to $\nu_{\Om,I}$.

On the figure, an intermediate result of this procedure is shown:
the set-theoretic difference of $K_0$ and a union of copies $\Om$
associated with one side of $K_0$.
 %\newpage
 
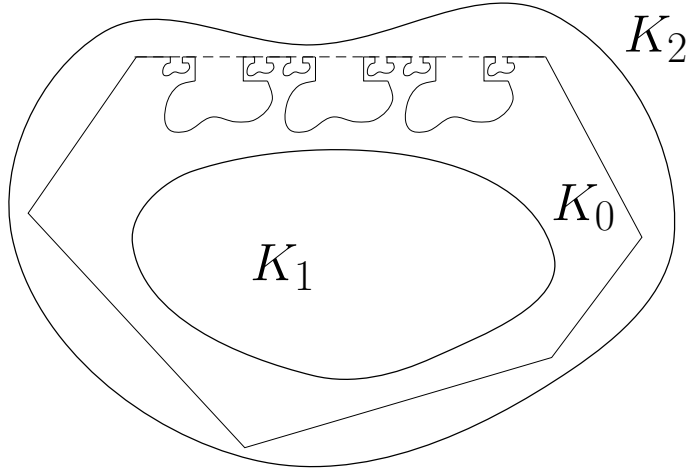
\vspace*{64mm}

       \scalebox{0.8}{
 \rput(-0.1,0){
\psline[linewidth=0.35pt](8.9,6)(10.5,3)(9,1)(3.9,-0.5)(0.3,3.4)(2.1,6)
\psline[linewidth=0.15pt,linestyle=dashed](8.9,6)(2.1,6)
                    }
                    \rput(7.9,6){
                    \scalebox{-0.4}{
              \rput(-2,0){
             \psline(1,1)(2,1)(2,0)
   \psline(4,0)(4,1)
 \pscurve(4,1)(4.7,1.2)(5.2,2)(5.1,3.0)(4.4,3.1)(3.3,2.5)(2.3,2.7)(1.5,2.6)(0.8,1.8)(1,1)
              }
             \rput(3,0){
             \psline(1,1)(2,1)(2,0)
   \psline(4,0)(4,1)
 \pscurve(4,1)(4.7,1.2)(5.2,2)(5.1,3.0)(4.4,3.1)(3.3,2.5)(2.3,2.7)(1.5,2.6)(0.8,1.8)(1,1)
              }
              \rput(8,0){
             \psline(1,1)(2,1)(2,0)
   \psline(4,0)(4,1)
 \pscurve(4,1)(4.7,1.2)(5.2,2)(5.1,3.0)(4.4,3.1)(3.3,2.5)(2.3,2.7)(1.5,2.6)(0.8,1.8)(1,1)
              }
    %%%%%%%%%%%%%%%%%%%%%%%%%
    %%%%%%%%%%%%%%%%%%%%%%%%%
     \rput(0.015,0){
           \rput(-1.5,0){
           \scalebox{0.25}{
 \psline[linewidth=3.2pt](1,1)(2,1)(2,0)
   \psline[linewidth=3.2pt](4,0)(4,1)
 \pscurve[linewidth=3.2pt](4,1)(4.7,1.2)(5.2,2)(5.1,3.0)(4.4,3.1)(3.3,2.5)(2.3,2.7)(1.5,2.6)(0.8,1.8)(1,1)
   }
   }
    \rput(2,0){
 \scalebox{0.25}{
  \psline[linewidth=3.2pt](1,1)(2,1)(2,0)
   \psline[linewidth=3.2pt](4,0)(4,1)
 \pscurve[linewidth=3.2pt](4,1)(4.7,1.2)(5.2,2)(5.1,3.0)(4.4,3.1)(3.3,2.5)(2.3,2.7)(1.5,2.6)(0.8,1.8)(1,1)
   }
   }
       \rput(3.5,0){
 \scalebox{0.25}{
  \psline[linewidth=3.2pt](1,1)(2,1)(2,0)
   \psline[linewidth=3.2pt](4,0)(4,1)
 \pscurve[linewidth=3.2pt](4,1)(4.7,1.2)(5.2,2)(5.1,3.0)(4.4,3.1)(3.3,2.5)(2.3,2.7)(1.5,2.6)(0.8,1.8)(1,1)
   }
   }
       \rput(7,0){
 \scalebox{0.25}{
  \psline[linewidth=3.2pt](1,1)(2,1)(2,0)
   \psline[linewidth=3.2pt](4,0)(4,1)
 \pscurve[linewidth=3.2pt](4,1)(4.7,1.2)(5.2,2)(5.1,3.0)(4.4,3.1)(3.3,2.5)(2.3,2.7)(1.5,2.6)(0.8,1.8)(1,1)
       }
       }
    \rput(8.5,0){
 \scalebox{0.25}{
  \psline[linewidth=3.2pt](1,1)(2,1)(2,0)
   \psline[linewidth=3.2pt](4,0)(4,1)
 \pscurve[linewidth=3.2pt](4,1)(4.7,1.2)(5.2,2)(5.1,3.0)(4.4,3.1)(3.3,2.5)(2.3,2.7)(1.5,2.6)(0.8,1.8)(1,1)
       }
       }
       \rput(12,0){
 \scalebox{0.25}{
  \psline[linewidth=3.2pt](1,1)(2,1)(2,0)
   \psline[linewidth=3.2pt](4,0)(4,1)
 \pscurve[linewidth=3.2pt](4,1)(4.7,1.2)(5.2,2)(5.1,3.0)(4.4,3.1)(3.3,2.5)(2.3,2.7)(1.5,2.6)(0.8,1.8)(1,1)
   }
   }
   \psline(-2.5,0)(-1,0)
   \psline(-0.5,0)(0,0)
   \psline(2,0)(2.5,0)
   \psline(3,0)(4,0)
   \psline(4.5,0)(5,0)
   \psline(7,0)(7.5,0)
   \psline(8,0)(9,0)
   \psline(9.5,0)(10,0)
   \psline(12,0)(12.5,0)
   \psline(13,0)(14.5,0)
       }
       }
       }
       \psecurve[linewidth=0.6pt](1.5,6.4)(5,6.2)(9,6.7)(11,3)(9.3,0.8)(4,-0.7)(0,3)(1.5,6.4)(5,6.2)(9,6.7)
        \psecurve[linewidth=0.6pt](3.1,4)(9,2.7)(7,1)(5,0.7)(2,2.9)(3,4.1)(9,2.7)(7,1)
        \rput(4.5,2.5){\Huge $K_1$}
         \rput(10.7,6.3){\Huge $K_2$}
           \rput(9.5,3.5){\Huge $K_0$}
                            }

\vspace{15mm}

The construction itself is simple, but its description is rather
cumbersome.  Let the polygon $K_0$ have $m$ sides; designate them by
$L_i$,\, $i = 1, \ldots, m$. On each segment $L_i$, put two points
at the distance $\ve/(4m)$ from its endpoints; the segment $L_i^\ve$
limited by these points has the length $|L_i| - \ve/(2m)$. Denote by
$\Pi_{i,\del,\ve}$ the closed rectangle that belongs to $K_0$, has a
side coincident with $L_i^\ve$, and another side of length $\del$
orthogonal to $L_i^\ve$. Choose $\del > 0$ in such a way that the
rectangles $\Pi_{i,\del,\ve}$,\, $i = 1, \ldots, m$ do not mutually
intersect and do not intersect with $K_1$.

Fix $i$ and let $\Pi^1 := \Pi_{i,\del,\ve}$. Select a ratio $k^1 >
0$ and a finite number of isometries $V_j^1$,\, $j = 1, \ldots, q_1$
in such a way that the sets $V_j^1 k^1 \tilde\Om^C$ do not mutually
intersect, the union of their closures $\cup_j \overline{V_j^1 k^1
\tilde\Om^C}$ is a rectangle containing $L_i^\ve$ and being
contained in $\Pi^1$, and all the segments $V_j^1 k^1 I$ belong to
$L_i^\ve$.

Denote $\Pi^2 := \cup_j V_j^1 k^1 \overline{(\tilde\Om_l \cup
\tilde\Om_r)}$. Select a ratio $k^2 > 0$ and a finite number of
isometries $V_j^2$,\, $j = 1, \ldots, q_2$ in such a way that the
sets $V_j^2 k^2 \tilde\Om^C$ do not mutually intersect, the set
$\cup_j \overline{V_j^2 k^2 \tilde\Om^C}$ is a union of rectangles
containing $L_i^\ve \cap \Pi^2$ and being contained in $\Pi^2$, and
all the segments $V_j^2 k^2 I$ belong to $L_i^\ve$.

Continuing this process, one obtains infinite sequences of sets
$\Pi^1,\ \Pi^2, \ldots$, isometries $\{ V_j^1 \},\ \{ V_j^2 \},
\ldots$, and numbers $k_1,\ k_2, \ldots$. The lengths
$$
 r_n := \Big|
L_i^\ve \setminus
 \left(
 \left( \cup_j V_j^1 k^1 I \right) \cup \ldots
 \cup \left( \cup_j V_j^n k^n I \right)
 \right)
 \Big|, \ \ \ n = 1,\, 2, \ldots
$$
form a decreasing geometric progression with ratio $1 - a/c$. Note
that if $a = c$, one has $r_1 = 0$, so the process terminates on the
first step.

Choose $n$ in such a way that $r_n \le \ve/(2m)$ and consider all
the sets $V_j^1 k^1 \Om$,\, $j = 1,\ldots, q_1$; \, $V_j^2 k^2
\Om$,\, $j = 1,\ldots, q_2$; \, $\ldots$; \, $V_j^n k^n \Om$,\, $j =
1,\ldots, q_n$, as well as the segments $V_j^1 k^1 I$,\, $j =
1,\ldots, q_1$; \, $V_j^2 k^2 I$,\, $j = 1,\ldots, q_2$; \,
$\ldots$; \, $V_j^n k^n I$,\, $j = 1,\ldots, q_n$. For further
convenience, rename these sets and segments according to
$\Om_{l,j}^i$ and $I_{l,j}^i,\ i = 1,\ldots, m,\ l = 1, \ldots, n,\
j = 1,\ldots, q_{li}$. The superscript $i$ enumerates the sides of
$K_0$, the subscript $l$ corresponds to the ''rank'' of the set on
the $i$th side, and the subscript $j$ enumerates identical sets of
rank $l$ on the $i$th side. The sets $\Om_{l,j}^i$ do not mutually
intersect and are copies of $\Om$, so the related measures are
identical, $\nu_{\Om_{l,j}^i,I_{l,j}^i} = \nu_{\Om,I}$. Besides,
each segment $I_{l,j}^i$ belongs to $\pl K_0$. Denote by
$|I_{l,j}^i|$ its length, denote $\kappa_{l,j}^i = |I_{l,j}^i|/|\pl
K_0|$ and let $\kappa_0 = 1 - \sum_{i,j,l} \kappa_{l,j}^i$ be the
part of the polygon's boundary free of these segments; one has
$\kappa_0 \le \ve/|\pl K_0|$.

Consider the set
$$
B_\ve = \overline{K_0 \setminus \left( \cup_{i,j,l} \Om_{l,j}^i
\right)};
$$
the bar means closure. One has $K_1 \subset B_\ve \subset K_2$,
hence $B_\ve \in \BBB_{K_1,K_2}$. Further, one has conv$B_\ve =
K_0$; the set $\pl(\text{conv}B_\ve) \setminus \pl B_\ve$ is a union
of a finite number of disjoint intervals $\cup_{i,j,l} I_{l,j}^i$.
One has
$$
\nu_{B_\ve} = \kappa_0 \nu^0 + \sum_{i,j,l} \kappa_{l,j}^i
\nu_{\Om_{l,j}^i,I_{l,j}^i}\,.
$$
Taking into account that $\sum_{i,j,l} \kappa_{l,j}^i = 1 -
\kappa_0$ and $\nu_{\Om_{l,j}^i,I_{l,j}^i} = \nu_{\Om,I}$, one
obtains that $\nu_{B_\ve} = \kappa_0 \nu^0 + (1 - \kappa_0)
\nu_{\Om,I}$. As $\ve \to 0$, one has $\kappa_0 \to 0$, therefore
$\nu_{B_\ve}$ converges in variation, and hence in the weak
topology, to $\nu_{\Om,I}$. This proves theorem \ref{t1}.

\section{Proof of the auxiliary result}

\subsection{A preparatory lemma}

Let $\s$ be a permutation of the set $\{ 1, \ldots, m \}$,\, $m \ge
2$. Denote $\theta_i^m = \arcsin (-1 + 2i/m)$,\, $i = 0,\, 1,
\ldots, m$. The points $\theta_1^m, \ldots, \theta_{m-1}^m$ divide
the segment $\interval$ into smaller segments $\Theta_i^m :=
[\theta_{i-1}^m,\, \theta_i^m]$ of equal measure $\lam$:\,\
$\lam(\Theta_i^m) = 2/m$,\, $i = 1, \ldots, m$. Recall that $\lam$
is defined by $d\lam(\vphi) = \cos\vphi\, d\vphi$. Define the
mapping $\vphi_\s:\, \interval \to \interval$ as follows: $\vphi_\s$
preserves the measure $\lam$ and rearranges the segments
$\Theta_i^m$,\, $i = 1, \ldots, m$, according to the permutation
$\s$; the restriction of $\vphi_\s$ to each segment $\Theta_i^m$
monotonically decreases and maps it bijectively (up to a set of zero
measure) to $\Theta_{\s(i)}^m$, if $\s(i) \ne i$, and is identity
mapping, if $\s(i) = i$. The mapping $\vphi_\s$ is given, up to a
set of zero measure, by the following formulas:
$$
\text{if} \ \vphi \in \Theta_i^m, \ \ \s(i) \ne i \ \ \text{then} \
\vphi_\s(\vphi) = \arcsin\left( -2 + 2 \frac{i + \s(i) - 1}{m} -
\sin\vphi \right); \ \ \ \
$$
$$
\text{if} \ \vphi \in \Theta_i^m, \ \ \s(i) = i \ \ \text{then} \
\vphi_\s(\vphi) = \vphi. \ \ \ \ \ \ \ \ \ \ \ \ \ \ \ \ \ \ \ \ \ \
\ \ \ \ \ \ \ \ \ \ \ \ \ \ \ \ \ \ \ \ \ \ \ \ \ \ \ \
$$

The mapping $\vphi_\s$ induces a Borel measure $\nu^\s$ on $Q$: for
any measurable set $A \subset Q$ one has
$$
\nu^\s(A) := \lam(\{ \vphi:\, (\vphi, \vphi_\s(\vphi)) \in A \}).
$$
The measure $\nu^\s$ satisfies the condition A1; if, additionally,
$\s^2 =$ id then it also satisfies A2. Thus, the set $\{ \nu^\s:\,
\s^2 = \text{id} \}$ belongs to $\MMM$. Later on, we will need the
following lemma.

\begin{lem}\label{dense}
The set $\{ \nu^\s:\, \s^2 = \mathrm{id},\, \s(1) \ne m \}$ is
everywhere dense in $\MMM$ in the weak topology.
\end{lem}

Its proof is not difficult and is put in appendix.

The condition $\s(1) \ne m$ is quite technical and will be used
below when proving corollary \ref{sledstvie c}.

\subsection{Definition of a reflector}\label{opredelenie reflektora}

Denote $e_\vphi = (\sin\vphi, \cos\vphi)$. Let $\vphi_1,\, \vphi_2
\in (-\pi/2,\, \pi/2)$,\, $\vphi_1 \ne \vphi_2$; consider two rays
$t e_{\vphi_1},\ t \ge 0$ and $t e_{\vphi_2},\ t \ge 0$. These rays
make angles $\vphi_1$ and $\vphi_2$ with the vector $e_0 = (0,1)$,
the angle being counted clockwise from $e_0$. Let $x^{(1)}$ and
$x^{(2)}$ be the points of intersection of these rays with the unit
circumference $|x - e_0| = 1$. Let $O = (0,0)$ denote the origin.
Draw two parabolas, $p_1$ and $p_2$, with common focus at $O$ and
with a common axis parallel to $x^{(2)} - x^{(1)}$; require,
moreover, that the parabola $p_1$ contains $x^{(1)}$, the parabola
$p_2$ contains $x^{(2)}$, and the intersection of convex sets
bounded by these parabolas contains the segment $[x^{(1)},\,
x^{(2)}]$.

Let $\del \ge 0$. Consider the straignt lines $l_0,\ l_+$ and $l_-$
given by the formulas $(x, e_0) = 0$,\, $(x, e_{-\del}) + \del
\sin\del = 0$ and $(x, e_{\del}) + \del \sin\del = 0$, respectively.
Let us denote by $R(\vphi_1, \vphi_2, \del)$ and call {\it
$(\vphi_1, \vphi_2, \del)$-reflector} the convex set bounded by
parabolas $p_1$,\, $p_2$ and lines $l_0$,\, $l_+$, and $l_-$, that
is
$$
R(\vphi_1, \vphi_2, \del) = \bigg\{ x:\ (x, e_0) \ge 0,\ (x,
e_{-\del}) + \del \sin\del \ge 0,\ (x, e_{\del}) + \del \sin\del \ge
0,
$$
\begin{equation}\label{o1}
\left. |x| - |x^{(1)}| \le \left(x - x^{(1)},\, \frac{x^{(2)} -
x^{(1)}}{|x^{(2)} - x^{(1)}|} \right),\ |x| - |x^{(2)}| \le \left(x
- x^{(2)},\, \frac{x^{(1)} - x^{(2)}}{|x^{(1)} - x^{(2)}|} \right)
\right\}.
\end{equation}

 
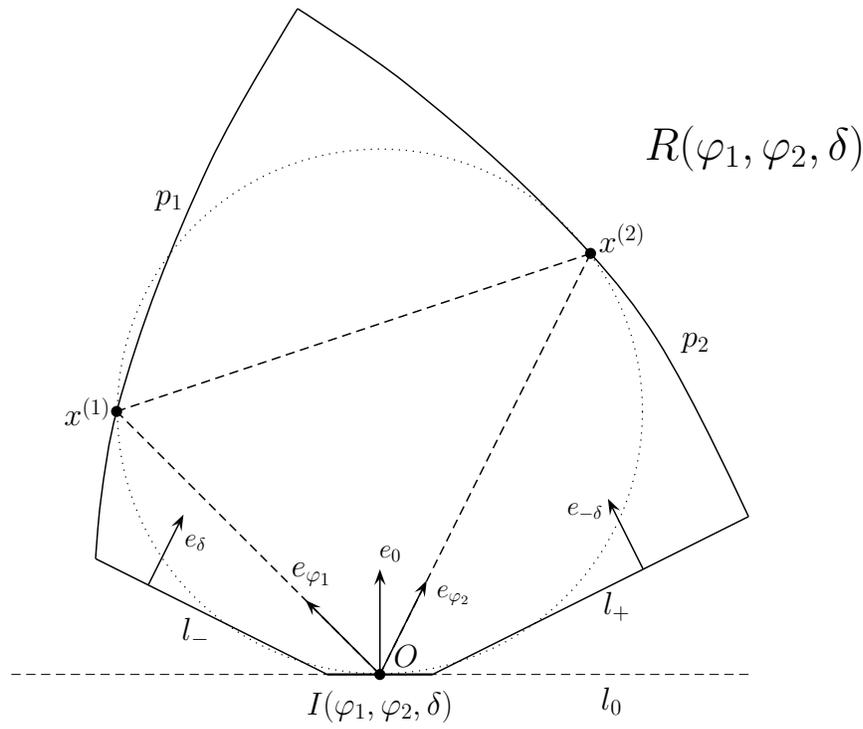
\vspace*{93mm}

 \scalebox{0.7}{
 \rput(8,0){
 \pscircle[linestyle=dotted](0,5){5}
 \psline[linestyle=dashed,linewidth=0.4pt](-7,0)(7,0)
  \psline[linewidth=1.2pt](-1,0)(1,0)
 \psline[linestyle=dashed](0,0)(-5,5)
 \psline[linestyle=dashed](0,0)(4,8)
 \psline[linestyle=dashed](-5,5)(4,8)
 \psline(-1,0)(-5.4,2.2)
 \psline(1,0)(7,3)
 \psdots[dotsize=6pt](-5,5)(4,8)(0,0)
 \pscurve(-5.4,2.2)(-5.3065,3.2704)(-5,5)(-3.2361,9.7082)(-1.5643,12.6547)
 \pscurve(7,3)(5.2489,6.3892)(4,8)(0.4633,11.2593)(-1.5643,12.6547)
 \rput(0,-0.6){\Large $I(\vphi_1,\vphi_2,\del)$}
  \rput(0.5,0.35){\Large $O$}
 \rput(-5.55,5){\Large $x^{(1)}$}
 \rput(4.6,8.3){\Large $x^{(2)}$}
       \rput(-4,9){\Large $p_1$}
       \rput(6,6.3){\Large $p_2$}
    \rput(4.4,-0.5){\Large $l_0$}
    \rput(-3.5,0.8){\Large $l_-$}
    \rput(4.5,1.3){\Large $l_+$}
     \psline[linewidth=0.8pt,arrows=->,arrowscale=2](0,0)(-1.414,1.414)
     \psline[linewidth=0.8pt,arrows=->,arrowscale=2](0,0)(0.89,1.788)
     \psline[linewidth=0.8pt,arrows=->,arrowscale=2](0,0)(0,2)
        \rput(-1.3,1.9){\Large $e_{\vphi_1}$}
        \rput(1.4,1.5){\large $e_{\vphi_2}$}
         \rput(0.2,2.3){\large $e_0$}
     \psline[linewidth=0.8pt,arrows=->,arrowscale=2](-4.4,1.7)(-3.73,3.04)
     \psline[linewidth=0.8pt,arrows=->,arrowscale=2](5,2)(4.33,3.34)
     \rput(3.9,3.1){\large $e_{-\del}$}
        \rput(-3.5,2.5){\large $e_\del$}
 }}
  \rput(10.5,7){\Large $R(\vphi_1, \vphi_2, \del)$}

 \vspace{17mm}

Let us call the set $I(\vphi_1, \vphi_2, \del) := \pl R(\vphi_1,
\vphi_2, \del) \cap l_0$ {\it base} of the reflector, and the set
$\kappa\, I(\vphi_1, \vphi_2, \del)$,\, $0 < \kappa < 1$, {\it
$\kappa$-base} of the reflector; for $\del$ small enough, these sets
coincide with $[-\del, \del] \times \{ 0 \}$ and $[-\kappa\del,
\kappa\del] \times \{ 0 \}$, respectively. The point $O$ is called
{\it center} of the reflector. In the particular case $\del = 0$,
the reflector $R(\vphi_1, \vphi_2) := R(\vphi_1, \vphi_2, 0)$ is the
set bounded by $p_1$,\, $p_2$, and $l_0$. Denote also
$I(\vphi_1,\vphi_2) := I(\vphi_1, \vphi_2, 0)$. Any copy of the
reflector, that is, any set of the form $V k\, R(\vphi_1, \vphi_2,
\del)$, where $V$ is a generic isometry and $k > 0$, will also be
called $(\vphi_1, \vphi_2, \del)$-reflector, the sets $V k\,
I(\vphi_1, \vphi_2, \del)$ and $V k\, \kappa\, I(\vphi_1, \vphi_2,
\del)$ will be called base and $\kappa$-base of this reflector,
respectively, and the point $V k\, O$ will be called center of the
reflector.

Note that $R(\vphi_1, \vphi_2, \del) = R(\vphi_2, \vphi_1, \del)$
and $R(\vphi_1, \vphi_2) = R(\vphi_2, \vphi_1)$.

Later on, we will need one more definition. Adding the left and
right hand sides of the inequalities in (\ref{o1}), one gets $2|x|
\le |x^{(1)}| + |x^{(2)}| + |x^{(2)} - x^{(1)}|$. The right hand
side of this inequality is perimeter of the triangle with vertices
$O$,\, $x^{(1)}$, and $x^{(2)}$, inscribed in the unit
circumference. It does not exceed the perimeter of an equilateral
triangle inscribed in the same circumference. This implies that
$2|x| \le 3\sqrt 3$. Define the set $\mathrm{T}(\del) = \{ x:\ 0 \le
(x, e_0) < 3\sqrt 3/2,\ (x, e_{-\del}) + \del \sin\del \ge 0,\ (x,
e_{\del}) + \del \sin\del \ge 0 \}$, the trapezium bounded from
below by $l_0$,\, $l_+$, and $l_-$, and from above, by a straight
line parallel to $l_0$. It follows from the construction that
$R(\vphi_1, \vphi_2, \del) \subset \mathrm{T}(\del)$ for any
$\vphi_1$ and $\vphi_2$.

 
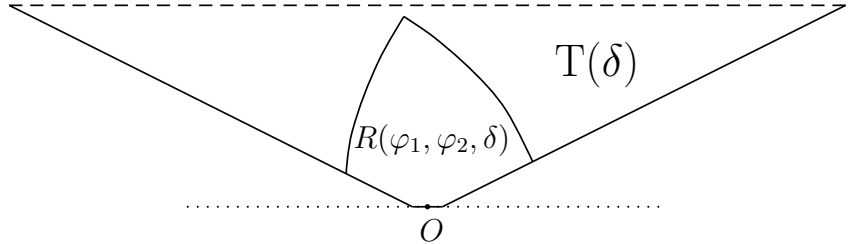
\vspace*{35mm}

 \scalebox{0.2}{
 \rput(35,0){
 \psline[linestyle=dotted,dotsep=15pt,linewidth=4pt](-16,0)(16,0)
  \psline[linewidth=3.2pt](-1,0)(1,0)
  \psline[linewidth=3.2pt](-1,0)(-27.8,13.4)
 \psline[linewidth=3.2pt](1,0)(27.8,13.4)
       \psline[linestyle=dashed,dash=25pt 15pt,linewidth=3.2pt](-27.8,13.4)(27.8,13.4)
 \psdots[dotsize=10pt](0,0)
 \pscurve[linewidth=3.2pt](-5.4,2.2)(-5.3065,3.2704)(-5,5)(-3.2361,9.7082)(-1.5643,12.6547)
 \pscurve[linewidth=3.2pt](7,3)(5.2489,6.3892)(4,8)(0.4633,11.2593)(-1.5643,12.6547)
     }}

 \rput(7.1,0.2){$O$}
  \rput(9.3,2.4){\Large $\mathrm{T}(\del)$}
   \rput(7.1,1.4){$R(\vphi_1, \vphi_2, \del)$}

 \vspace{7mm}

Consider the billiard in $R(\vphi_1,\vphi_2)$. Suppose that a
billiard particle is initially located at a regular point of the
boundary, $x = (\xi, 0) \in I(\vphi_1,\vphi_2)$, and has velocity
$e_{\vphi}$,\, $-\pi/2 < \vphi < \pi/2$, and during the subsequent
motion reflects first from $p_1$, then from $p_2$, and for the third
time, from an interior point of $I(\vphi_1,\vphi_2)$. Denote by $x^+
= (\xi^+, 0)$ the point of third reflection and denote by
$-e_{\vphi^+}$,\, $-\pi/2 < \vphi^+ < \pi/2$ the particle's velocity
just before the third reflection. The set of values $(\xi, \vphi,
\vphi_1, \vphi_2)$ such that the mentioned order of reflections
takes place is denoted by $\AAA$. The set $\AAA \subset \RRR^4$ is
open and contains the points of kind $(0, \vphi_1, \vphi_1,
\vphi_2)$,\, $\vphi_1 \ne \vphi_2$; hence it is nonempty. Denote by
$\tilde \xi^+$,\, $\tilde \vphi^+$ the mappings that send $(\xi,
\vphi, \vphi_1, \vphi_2) \in \AAA$ to $\xi^+$,\, $\vphi^+$,
respectively, and define the mapping $\tilde T:\, \AAA \to \RRR^2$
by $\tilde T(\xi, \vphi, \vphi_1, \vphi_2) = (\tilde \xi^+(\xi,
\vphi, \vphi_1, \vphi_2),\, \tilde \vphi^+(\xi, \vphi, \vphi_1,
\vphi_2))$. Recalling the definition of functions $x^+_{\Om,I}$,\,
$\vphi^+_{\Om,I}$ given in section \ref{auxiliary}, one concludes
that $(\tilde \xi^+(\xi, \vphi, \vphi_1, \vphi_2),\, 0) =
x^+_{R(\vphi_1, \vphi_2),I(\vphi_1, \vphi_2)}(x,\vphi)$,\, \ $\tilde
\vphi^+(\xi, \vphi, \vphi_1, \vphi_2) = \vphi^+_{R(\vphi_1,
\vphi_2),I(\vphi_1, \vphi_2)}(x,\vphi)$,\, where $x = (\xi, 0)$.

The mapping $\tilde T$ is infinitely differentiable; besides, one
has
$$
\tilde\xi^+(0, \vphi, \vphi_1, \vphi_2) = 0, \ \ \ \ \
\tilde\vphi^+(0, \vphi_1, \vphi_1, \vphi_2) = \vphi_2.
$$
Indeed, a particle gets out of the common focus $O$ of the parabolas
$p_1$ and $p_2$, then reflects from $p_1$ and moves in parallel to
the common axis of $p_1$ and $p_2$; finally, it reflects from $p_2$
and gets into the common focus $O$. If, besides, the initial
velocity equals $e_{\vphi_1}$ then after reflections at the points
$x^{(1)}$ and $x^{(2)}$ it returns to $O$ with velocity
$-e_{\vphi_2}$.

\begin{lem}\label{l3}
For any $\vphi_1,\, \vphi_2 \in (-\pi/2,\, \pi/2)$,\, $\vphi_1 \ne
\vphi_2$ holds
$$
\text{(a)} \ \ \ \ \ \ \
\frac{\pl\tilde\vphi^+}{\pl\vphi}\Bigg\rfloor _{\!\!\!\scriptsize
\begin{array}{l}\xi=0\\ \text{\raisebox{1mm}{$\vphi=\vphi_1$}} \end{array}}
 \!\!\!\!\!\!
(\xi,\vphi,\vphi_1,\vphi_2)\, =\,
-\frac{\cos\vphi_1}{\cos\vphi_2}\,; \ \ \ \ \ \ \ \ \
$$
$$
\text{(b)} \ \ \ \ \ \ \ \Big|\,
\frac{\pl\tilde\xi^+}{\pl\xi}\Bigg\rfloor _{\!\!\!\scriptsize
\begin{array}{l}\xi=0\\ \text{\raisebox{1mm}{$\vphi=\vphi_1$}} \end{array}}
 \!\!\!\!\!\!
(\xi,\vphi,\vphi_1,\vphi_2)\, \Big| =\, 1. \ \ \ \ \ \ \ \ \ \ \ \ \
\ \ \
$$
\end{lem}

\begin{proof}
Fix the values $-\pi/2 < \vphi_1 < \vphi_2 < \pi/2$; the case
$\vphi_1 > \vphi_2$ is completely analogous. Put $A = x^{(1)}$,\, $B
= x^{(2)}$. Note that on the figure below there are drawn the angles
$\vphi_1 < 0$,\, $\vphi_2 > 0$.

A particle that gets out of $O$ with velocity $e_{\vphi_1 +
\Del\vphi}$, with $\Del\vphi$ sufficiently small, first reflects
from $p_1$ at a point $\check{x}^{(1)} = A'$, then reflects from
$p_2$ at a point $\check{x}^{(2)} = B'$, and finally returns to $O$,
the velocity just before returning being $e_{\vphi_2 +
\Del\vphi^+}$. Thus, one has $\vphi_2 + \Del\vphi^+ =
\tilde\vphi^+(0,\vphi_1 + \Del\vphi,\vphi_1,\vphi_2)$. Put
$\Del\vphi < 0$, then $\Del\vphi^+ > 0$. On the figure, there is
shown the circumference $|x - e_0| = 1$. Note that the points $A$
and $B$ belong to the circumference, but $A'$ and $B'$, generally,
do not.

For our convenience, introduce two additional points $L_1,\ L_2$
lying on the the line $l_0$, one to the left and another to the
right of $O$. On can put, for example, $L_1 = (-\xi_0,0)$,\, $L_2 =
(\xi_0,0)$, where $\xi_0 > 0$. One has
$$
\measuredangle AOB = \vphi_2 - \vphi_1;
$$
$\measuredangle BOL_2 = \pi/2 - \vphi_2$; moreover, the angles
$\measuredangle OAB$ and $\measuredangle BOL_2$ are inscribed in the
same arc and hence are equal, therefore
$$
\measuredangle OAB = \pi/2 - \vphi_2.
$$
Further, $\measuredangle AOL_1 = \vphi_1 + \pi/2$; the angles
$\measuredangle ABO$ and $\measuredangle AOL_1$ are inscribed in the
same arc and hence are equal; therefore,
$$
\measuredangle ABO = \vphi_1 + \pi/2.
$$

 %\newpage

 
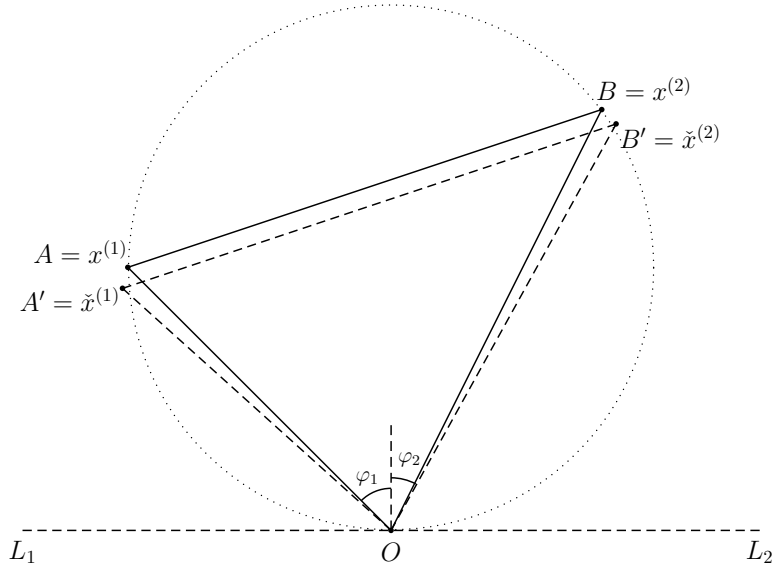
\vspace*{76mm}

 \scalebox{0.7}{
 \rput(7,0){
 \pscircle[linestyle=dotted](0,5){5}
 \psline[linestyle=dashed](-7,0)(7,0)
 \psline(0,0)(-5,5)
 \psline(0,0)(4,8)
 \psline(-5,5)(4,8)
        \psline[linestyle=dashed](0,0)(-5.1,4.6)
        \psline[linestyle=dashed](-5.1,4.6)(4.275,7.725)
        \psline[linestyle=dashed](0,0)(4.275,7.725)
        \psdots[dotsize=3pt](-5,5)(4,8)(0,0)(-5.1,4.6)(4.275,7.725)
 \rput(0,-0.4){\large $O$}
    \rput(-5.9,5.3){\large $A = x^{(1)}$}
 \rput(4.8,8.4){\large $B = x^{(2)}$}
    \rput(-6.1,4.4){\large $A' = \check{x}^{(1)}$}
 \rput(5.3,7.5){\large $B' = \check{x}^{(2)}$}
     \psline[linewidth=0.8pt,linestyle=dashed](0,0)(0,2)
       \psarc(0,0){0.8}{90}{135}
        \psarc(0,0){1.0}{62}{90}
        \rput(-0.45,1.0){$\vphi_1$}
        \rput(0.35,1.35){$\vphi_2$}
  \rput(-7,-0.4){\large $L_1$}
  \rput(7,-0.4){\large $L_2$}
 }}

 \vspace{20mm}

The law of sines implies that
\begin{equation}\label{trig1}
\frac{|OA|}{|OB|} = \frac{\sin\measuredangle ABO}{\sin\measuredangle
OAB} = \frac{\cos\vphi_1}{\cos\vphi_2}\,.
\end{equation}

By the elastic reflection law, the inner normal to $\pl
R(\vphi_1,\vphi_2)$ at $A$ makes the angle $\frac 12 \measuredangle
OAB = \pi/4 - \vphi_2/2$ with the straight line  $OA$; hence the
tangent to $\pl R(\vphi_1,\vphi_2)$ at $A$ makes the angle $\pi/4 +
\vphi_2/2$ with  $OA$. This implies that
$$
\measuredangle OAA' = \frac{\pi}{4} + \frac{\vphi_2}{2} +
O(\Del\vphi);
$$
besides
$$
\measuredangle OA'A = \pi - \measuredangle OAA' - \measuredangle
AOA' = \frac{3\pi}{4} - \frac{\vphi_2}{2} + O(\Del\vphi).
$$
(Note that $\measuredangle AOA' = -\Del\vphi$). Applying the law of
sines to the triangle $OAA'$, one gets
$$
\frac{\sin \measuredangle OA'A}{|OA|} = \frac{\sin \measuredangle
AOA'}{|AA'|}\,,
$$
hence
\begin{equation}\label{trig2}
\frac{\sin \left( 3\pi/4 - \vphi_2/2 + O(\Del\vphi) \right)}{|OA|} =
\frac{-\Del\vphi \left( 1 + O(\Del\vphi) \right)}{|AA'|}\,, \ \ \ \
\Del\vphi \to 0.
\end{equation}

Further, two different particles starting at the focus $O$, after
reflecting off the parabola $p_1$ move in parallel, that is, the
straight lines $AB$ and $A'B'$ are parallel. Denote the distance
between them by $\Del x$. One has
$$ \Del x = |AA'| \cdot \sin
\measuredangle A'AB = |AA'| \cdot \sin \left( \measuredangle OAA' +
\measuredangle OAB \right) =
$$
$$
= |AA'| \cdot \sin \left( 3\pi/4 - \vphi_2/2 + O(\Del\vphi) \right),
\ \ \ \ \ \ \Del\vphi \to 0.
$$
Taking account of (\ref{trig2}), one obtains
\begin{equation}\label{trig3}
\Del x = -|OA| \cdot \Del\vphi \cdot \left( 1 + O(\Del\vphi)
\right), \ \ \ \ \ \Del\vphi \to 0.
\end{equation}
Next, considering the triangle $OBB'$ in an analogous way and taking
into account that $\Del\vphi^+ > 0$, one obtains
\begin{equation}\label{trig4}
\Del x = |OB| \cdot \Del\vphi^+ \cdot \left( 1 + O(\Del\vphi)
\right), \ \ \ \ \ \Del\vphi \to 0.
\end{equation}
From (\ref{trig3}) and (\ref{trig4}) one gets
$$
\frac{\Del\vphi^+}{\Del\vphi} = -\frac{|OA|}{|OB|} \cdot \left( 1 +
O(\Del\vphi) \right), \ \ \ \ \Del\vphi \to 0;
$$
taking the limit in this equality as $\Del\vphi \to 0$ and using
(\ref{trig1}), one obtains
\begin{equation}\label{trig5}
\frac{\pl\tilde\vphi^+}{\pl\vphi} (\xi = 0,\vphi =
\vphi_1,\vphi_1,\vphi_2) = -\frac{\cos\vphi_1}{\cos\vphi_2}\,.
\end{equation}
The statement (a) of lemma is proved.

Further, one has $\tilde\xi^+(0,\vphi,\vphi_1,\vphi_2) = 0$; hence
for any $\vphi$ one has
\begin{equation}\label{trig6}
\frac{\pl\tilde\xi^+}{\pl\vphi}(0,\vphi,\vphi_1,\vphi_2) = 0.
\end{equation}
The restriction of $\tilde T$ to the subspace $\{ \vphi_1 =
\text{const},\, \vphi_2 = \text{const} \}$ preserves the measure
$d\mu(\xi,\vphi) = \cos\vphi\, d\xi d\vphi$, that is,
$\cos\tilde\vphi^+\, d\tilde\xi^+ d\tilde\vphi^+ = \cos\vphi\, d\xi
d\vphi$, and using that $d\tilde\xi^+ d\tilde\vphi^+ = \big|
\frac{D(\tilde\xi^+,\tilde\vphi^+)}{D(\xi,\vphi)} \big|\, d\xi
d\vphi$, one obtains
\begin{equation}\label{trig7}
\cos\tilde\vphi^+(\xi,\vphi,\vphi_1,\vphi_2) \cdot \bigg|
\frac{D(\tilde\xi^+,\tilde\vphi^+)}{D(\xi,\vphi)} \bigg| =
\cos\vphi.
\end{equation}
Taking into account that
$$
\frac{D(\tilde\xi^+,\tilde\vphi^+)}{D(\xi,\vphi)} =
\frac{\pl\tilde\xi^+}{\pl\xi} \ \frac{\pl\tilde\vphi^+}{\pl\vphi} -
\frac{\pl\tilde\xi^+}{\pl\vphi} \ \frac{\pl\tilde\vphi^+}{\pl\xi}
$$
and using (\ref{trig5}) and (\ref{trig6}), one gets
$$
\bigg| \frac{D(\tilde\xi^+,\tilde\vphi^+)}{D(\xi,\vphi)}
\bigg|_{(\xi=0,\vphi=\vphi_1,\vphi_1,\vphi_2)} = \Big|
\frac{\pl\tilde\xi^+}{\pl\xi}(0,\vphi_1,\vphi_1,\vphi_2) \Big| \cdot
\frac{\cos\vphi_1}{\cos\vphi_2}\,.
$$
Substituting this value into (\ref{trig7}) and using that $\tilde
\vphi^+(0, \vphi_1, \vphi_1, \vphi_2) = \vphi_2$, one obtains
\begin{equation}\label{trig8}
\Big| \frac{\pl\tilde\xi^+}{\pl\xi}(0, \vphi_1, \vphi_1, \vphi_2)
\Big| = 1;
\end{equation}
thus, the statement (b) of lemma is also proved. Using an additional
geometric argument, the formula (\ref{trig8}) could be made more
precise: $\,
\frac{\pl\tilde\xi^+}{\pl\xi}(\xi=0,\vphi=\vphi_1,\vphi_1,\vphi_2) =
1$, but we will not need this specification in the future.
\end{proof}

Let $f$ be a function of $x,\, \vphi,\, \vphi_1,\, \vphi_2,\, \del$.
In the following lemma, by $O(f)$ we mean a generic function of the
same variables such that $O(f)/f$ is bounded over all admissible
values of the variables.

Let $-\pi/2 < \Phi_1 < \Phi_2 < \pi/2$,\, $0 < \Phi_0 < \Phi_2 -
\Phi_1$.

\begin{lem}\label{reflector lemma}
There exist positive constants $c = c(\Phi_0, \Phi_1, \Phi_2)$ and
$c_0 = c_0(\Phi_0, \Phi_1, \Phi_2)$ such that
 \vspace{1mm}

(a)\, if $\vphi_1,\, \vphi_2 \in [\Phi_1,\, \Phi_2]$,\, $|\vphi_2 -
\vphi_1| \ge \Phi_0$,\, $|x| \le \del (1 - c\del -c|\vphi -
\vphi_1|)$ then
\begin{equation}\label{ref lemma phi}
\vphi^+_{R(\vphi_1,\vphi_2,\del),I(\vphi_1,\vphi_2,\del)} (x, \vphi)
- \vphi_2 = -\frac{\cos\vphi_1}{\cos\vphi_2}\, (\vphi - \vphi_1) +
O(\del) + O((\vphi - \vphi_1)^2);
\end{equation}

(b)\, if, moreover, $\del < c_0$,\, $|\vphi - \vphi_1| < c_0$ then
\begin{equation}\label{ref lemma n}
n_{R(\vphi_1,\vphi_2,\del),I(\vphi_1,\vphi_2,\del)} (x, \vphi) = 3.
\end{equation}
\end{lem}

Thus, the lemma states that a billiard particle that starts moving
at a point of a reflector's base, under the given restrictions on
the initial position $x$, on the initial direction of motion
$\vphi$, and on the reflector's parameters, after two reflections
from $p_1$ and $p_2$ will return to the reflector's base and the
direction of motion of the returning particle $\vphi^+$ will satisfy
the relation (\ref{ref lemma phi}).

\begin{proof}
The set $\{ (0, \vphi_1, \vphi_1, \vphi_2):\, \vphi_1,\, \vphi_2 \in
[\Phi_1,\, \Phi_2],\ |\vphi_2 - \vphi_1| \ge \Phi_0 \}$ is compact
and belongs to the open set $\AAA$; hence there exist positive
values $\Del\xi$ and $\Del\vphi$ such that the compact set
$\AAA(\Del\xi, \Del\vphi, \Phi_0, \Phi_1, \Phi_2) := \{ (\xi, \vphi,
\vphi_1, \vphi_2):\, |\xi| \le \Del\xi,\, |\vphi - \vphi_1| \le
\Del\vphi;\ \vphi_1,\, \vphi_2 \in [\Phi_1,\, \Phi_2],\ |\vphi_2 -
\vphi_1| \ge \Phi_0 \}$ also belongs to $\AAA$.

For $(\xi, \vphi, \vphi_1, \vphi_2) \in \AAA(\Del\xi, \Del\vphi,
\Phi_0, \Phi_1, \Phi_2)$ holds
$$
\tilde\vphi^+(\xi, \vphi, \vphi_1, \vphi_2)\, -\, \tilde\vphi^+(0,
\vphi_1, \vphi_1, \vphi_2) =
$$
$$
= [\tilde\vphi^+(\xi, \vphi, \vphi_1, \vphi_2) - \tilde\vphi^+(0,
\vphi, \vphi_1, \vphi_2)]\ +\ [\tilde\vphi^+(0, \vphi, \vphi_1,
\vphi_2) - \tilde\vphi^+(0, \vphi_1, \vphi_1, \vphi_2)] =
$$
$$
= \frac{\pl\tilde\vphi^+}{\pl\xi} \cdot \xi\, +\,
\frac{\pl\tilde\vphi^+}{\pl\vphi}(0, \vphi_1, \vphi_1, \vphi_2)
\cdot (\vphi - \vphi_1)\, +\, \frac{\pl^2 \tilde\vphi^+}{\pl\vphi^2}
\cdot \frac{(\vphi - \vphi_1)^2}{2}\,.
$$
According to statement (a) of lemma \ref{l3}, one has
$\frac{\pl\tilde\vphi^+}{\pl\vphi} (0, \vphi_1, \vphi_1, \vphi_2) =
-\frac{\cos\vphi_1}{\cos\vphi_2}$. Moreover, using that
$\tilde\vphi^+(0, \vphi_1, \vphi_1, \vphi_2) = \vphi_2$ and that the
functions $\frac{\pl\tilde\vphi^+}{\pl\xi}$,\, $\frac{\pl^2
\tilde\vphi^+}{\pl\vphi^2}$ are bounded on the set $\AAA(\Del\xi,
\Del\vphi, \Phi_0, \Phi_1, \Phi_2)$, and introducing the shorthand
notation $\tilde\vphi^+(\xi, \vphi, \vphi_1, \vphi_2) =
\tilde\vphi^+$, one obtains
\begin{equation}\label{refl1}
\tilde\vphi^+ - \vphi_2 = -\frac{\cos\vphi_1}{\cos\vphi_2}\, (\vphi
- \vphi_1) + O(\xi) + O \left( (\vphi - \vphi_1)^2 \right).
\end{equation}

Further, one has
$$
\tilde\xi^+(\xi, \vphi, \vphi_1, \vphi_2)\, -\, \tilde\xi^+(0,
\vphi, \vphi_1, \vphi_2)\, =\, \frac{\pl\tilde\xi^+}{\pl\xi} (0,
\vphi, \vphi_1, \vphi_2) \cdot \xi\, +\, \frac{\pl^2
\tilde\xi^+}{\pl\xi^2} \cdot \frac{\xi^2}{2}\,,
$$
$$
\frac{\pl\tilde\xi^+}{\pl\xi} (0, \vphi, \vphi_1, \vphi_2)\, =\,
\frac{\pl\tilde\xi^+}{\pl\xi} (0, \vphi_1, \vphi_1, \vphi_2)\, +\,
\frac{\pl^2 \tilde\xi^+}{\pl\xi \pl\vphi} \cdot (\vphi - \vphi_1).
$$
According to statement (b) of lemma \ref{l3},\
$|\frac{\pl\tilde\xi^+}{\pl\xi} (0, \vphi_1, \vphi_1, \vphi_2)| =
1$. Taking into account that $\tilde\xi^+(0, \vphi, \vphi_1,
\vphi_2) = 0$ and the values $\frac{\pl^2 \tilde\xi^+}{\pl\xi^2}$,\,
$\frac{\pl^2 \tilde\xi^+}{\pl\xi \pl\vphi}$ are bounded on
$\AAA(\Del\xi, \Del\vphi, \Phi_0, \Phi_1, \Phi_2)$, and using the
shorthand notation $\tilde\xi^+(\xi, \vphi, \vphi_1, \vphi_2) =
\tilde\xi^+$, one gets
\begin{equation}\label{refl2}
\tilde\xi^+ = \xi \left( \pm 1 + O(\xi) + O(\vphi - \vphi_1)
\right).
\end{equation}

It follows from (\ref{refl1}) that for some constant $c_0 > 0$ the
following holds: if $|\xi|$,\, $|\vphi - \vphi_1| < c_0$ then
$|\vphi|$,\, $|\tilde\vphi^+| < \pi/2 - c_0$. From (\ref{refl2}) it
follows that for some constant $c > 0$ the following holds: if
$|\xi| \le \del (1 - c\del -c|\vphi - \vphi_1|)$ then $|\tilde\xi^+|
< \del$. These facts imply that for $\del < c_0$ the following holds
true: if $|\xi| \le \del (1 - c\del -c|\vphi - \vphi_1|)$ and
$|\vphi - \vphi_1| < c_0$ then $|\tilde\xi^+| < \del$ and
$|\vphi|$,\, $|\tilde\vphi^+| < \pi/2 - c_0$.

The last conclusion means that under the chosen restrictions on the
initial data $\xi$,\, $\vphi$, a billiard particle in
$R(\vphi_1,\vphi_2)$ starts moving at some point on the base
$I(\vphi_1,\vphi_2)$ and after two reflections off $p_1$ and $p_2$
returns to the base, besides the distances from the initial and
final points to the origin $O$ are less than $\del$ and the
directions of initial and final motion make with the base angles
more than $\del$. This implies that this motion is at the same time
a billiard motion in $R(\vphi_1,\vphi_2,\del)$, that is,
$(\tilde\xi^+, 0) =
x^+_{R(\vphi_1,\vphi_2,\del),I(\vphi_1,\vphi_2,\del)} (x, \vphi)$,\,
$\tilde\vphi^+ = \vphi^+_{R(\vphi_1, \vphi_2, \del),I(\vphi_1,
\vphi_2, \del)} (x, \vphi)$, and $n_{R(\vphi_1, \vphi_2, \del),
I(\vphi_1, \vphi_2, \del)} (x, \vphi) = 3$ (recall that $x = (\xi,
0)$). Thus, (\ref{ref lemma n}) is proved.

Taking into account the above-mentioned and using (\ref{refl1}), one
concludes that (\ref{ref lemma phi}) is true locally, for $\del <
c_0$,\, $|\vphi - \vphi_1| < c_0$; hence it is also true globally,
for all admissible values of $x,\, \vphi,\, \vphi_1,\, \vphi_2,\,
\del$.
\end{proof}

\subsection{Definition of a pair $(\Om,I)$}\label{opredelenie pary}

Consider a permutation $\s$ of $\{ 1, \ldots, m \}$ such that $\s^2
= $id and $\s(1) \ne m$. Fix $r > 1$ and $\del > 0$, and put $I =
[-1/2,\, 1/2] \times \{ 0 \}$. In this section we shall define a set
$\Om = \Om(\s, r, \del)$ such that $\pl\Om$ contains $I$, and in the
next section we shall prove that $\nu_{\Om(\s, r, \del),I}$ weakly
converges to $\nu^\s$ as $r \to \infty$,\, $r\del \to 0$.

Recall that by definition, $\theta_i = \theta_i^m = \arcsin (-1 +
2i/m)$,\, $\Theta_i = \Theta_i^m = [\theta_{i-1},\, \theta_i]$,\, $i
= 0,\ 1, \ldots, m$, and $e_\vphi = (\sin\vphi, \cos\vphi)$. Denote
$P_i = r e_{\theta_i}$; the points $P_i$ lie on the upper
half-circumference of radius $r$ with the center $O = (0,0)$. Denote
by $\QQQ$ the polygon $P_0 P_1 \ldots P_m P_0$. Denote by $n_i$ the
unit outer normal to $\pl\QQQ$ at some point of $[P_{i-1},\, P_i]$;
one has $n_i = e_{(\theta_{i-1}+\theta_i)/2}$.

 
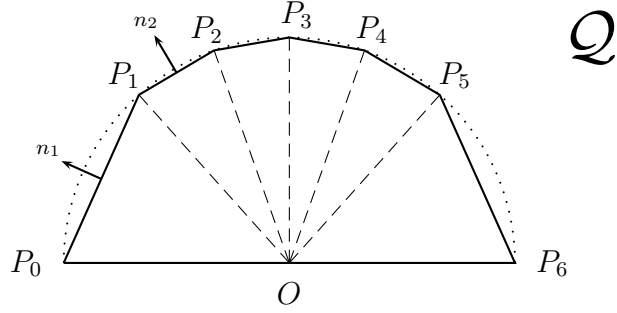
\vspace*{36mm}

 \rput(6,-0.5){
 \psarc[linestyle=dotted,linewidth=0.8pt](0,0){3}{0}{180}
 \psline(-3,0)(-2,2.236)(-1,2.828)(0,3)(1,2.828)(2,2.236)(3,0)(-3,0)
               \psline{->}(-2.5,1.118)(-3.031,1.3555)
                 \rput(-3.2,1.5){\scriptsize $n_1$}
               \psline{->}(-1.5,2.532)(-1.796,3.032)
                \rput(-1.9,3.2){\scriptsize $n_2$}
 \psline[linestyle=dashed,linewidth=0.4pt](0,0)(-2,2.236)
  \psline[linestyle=dashed,linewidth=0.4pt](0,0)(-1,2.828)
   \psline[linestyle=dashed,linewidth=0.4pt](0,0)(0,3)
    \psline[linestyle=dashed,linewidth=0.4pt](0,0)(1,2.828)
     \psline[linestyle=dashed,linewidth=0.4pt](0,0)(2,2.236)
  \rput(0,-0.4){$O$}
  \rput(-3.5,0){$P_0$}
  \rput(-2.2,2.5){$P_1$}
  \rput(-1.1,3.1){$P_2$}
  \rput(3.5,0){$P_6$}
  \rput(2.2,2.5){$P_5$}
  \rput(1.1,3.1){$P_4$}
  \rput(0.1,3.3){$P_3$}
       \rput(4,3){\Huge $\QQQ$}
  }
 \vspace{17mm}

Fix the values $i,\, j \in \{ 1, \ldots, m \}$,\, $i \ne j$, and a
point $x \in [P_{i-1},\, P_i]$. Denote by $\theta$ the angle the
vector $x$ makes with $e_0$; thus, $x = |x| e_\theta$,\, $\theta \in
\Theta_i$. Denote by $\vphi_{ij}$ the bijective mapping from
$\Theta_i$ to $\Theta_j$ that monotonically decreases and preserves
the measure $\lam$. Put $\theta' = \vphi_{ij}(\theta)$; one easily
sees that $\sin\theta_j - \sin\theta' = \sin\theta -
\sin\theta_{i-1}$. Take the vector $x' = x'(x) \in [P_{j-1},\, P_j]$
that forms the angle $\theta'$ with $e_0$; thus, $x' = |x'|
e_{\theta'}$. Denote by $\vphi_1 = \vphi_1(x)$ and $\vphi_2 =
\vphi_2(x)$ the angles the vectors $x$ and $x - x'$, respectively,
form with $n_i$; these angles are counted from $n_i$ clockwise; in
particular, one has $\vphi_1 = \theta - (\theta_{i-1}+\theta_i)/2$.
Finally, define the reflector $R_{ij}^k[x] = R_{ij}^k[x,r,\del]$ by
$$
R_{ij}^k[x] = V_x k R\left( \vphi_1, \vphi_2, \del \right),
$$
where $0 < k < 1$ and $V_x$ is the isometry taking the point $O$ to
the point $x$, and the vector $e_0$, to the vector $n_i$. Note that
for fixed $m$ and under the condition that $x$ belongs to the broken
line $P_0 P_1 \ldots P_m$, the value $x$ uniquely defines $r$, and
hence the isometry $V_x$.

Let us give an example explaining the meaning of the above
definition. Consider the billiard in the set $\QQQ \cup
R_{ij}^{k_1}[x] \cup R_{ji}^{k_2}[x']$: the union of the polygon and
two reflectors; the ratios $k_1$ and $k_2$ being taken small enough,
so that the reflectors do not intersect. If a particle starts at $O$
with the initial velocity $e_\theta$, then in the subsequent motion
it passes through $x$, makes two reflections off the boundary of the
first reflector, then moves along the straight line joining the
points $x$ and $x'$, makes two more reflections off the boundary of
the second reflector, and finally returns to $O$, the velocity just
before the return being $-e_{\theta'}$ (see the figure below).
 %\newpage

 
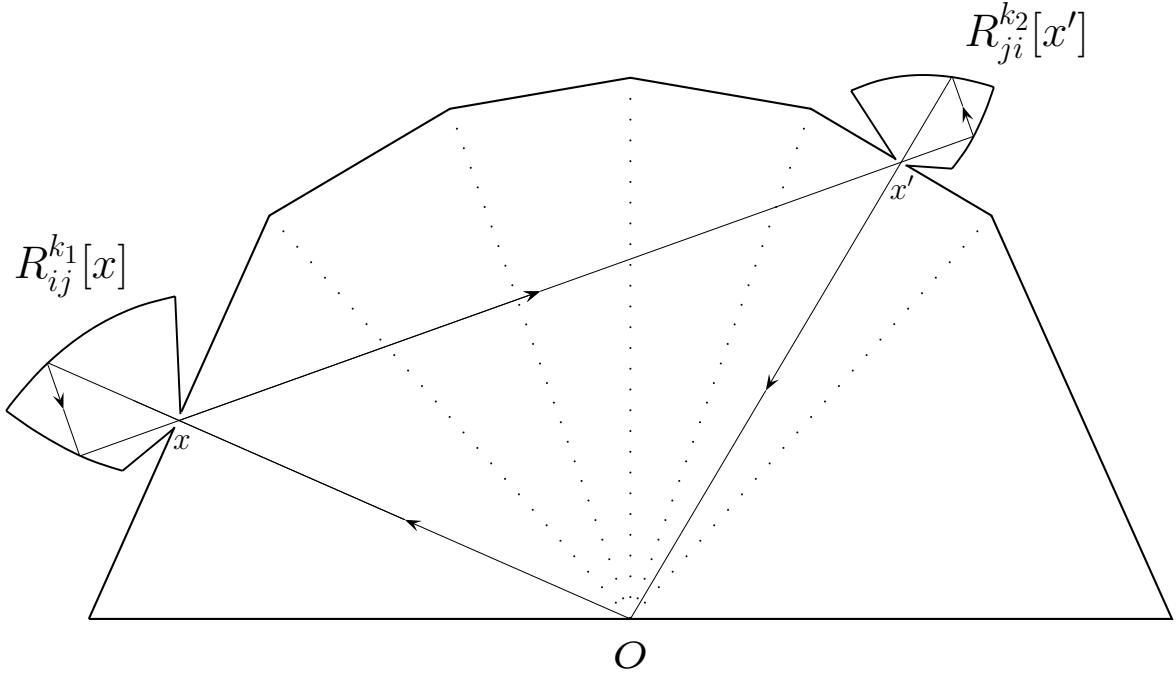
\vspace*{80mm}

\rput(-0.8,0){
  \scalebox{2.4}{
\rput(3,-0.5){
 \psline[linewidth=0.16pt](0,0)(0.75,1.266)
 \psline[linewidth=0.16pt,arrows=->,arrowscale=1](0,0)(-1.25,0.55)
 \psline[linewidth=0.16pt](-1.25,0.55)(-2.75,1.21)
         \psline[linewidth=0.16pt](-1.25,0.55)(-3.23125,1.42175)
 \psline[linewidth=0.16pt,arrows=<-,arrowscale=1](0.75,1.266)(1.782,3.008016)%(1.65,2.7852)
 \psline[linewidth=0.16pt,arrows=->,arrowscale=1](-2.5,1.1)(-0.5,1.816)
        \psline[linewidth=0.16pt,arrows=->,arrowscale=1](-3.05,0.9031)(-0.5,1.816)
        \psline[linewidth=0.16pt,arrows=->,arrowscale=1](-3.23125,1.42175)(-3.140625,1.162425)
        \psline[linewidth=0.16pt](-3.140625,1.162425)(-3.05,0.9031)
\psline[linewidth=0.16pt](-0.5,1.816)(1.9,2.6752)%(1.5,2.532)
                     \psline[linewidth=0.16pt,arrows=->,arrowscale=1](1.9,2.6752)(1.841,2.841608)
                     \psline[linewidth=0.16pt](1.841,2.841608)(1.782,3.008016)
\psline[linewidth=0.33pt](-3,0)(-2.525,1.0641)
  \psline[linewidth=0.33pt](-2.4875,1.148)(-2,2.236)(-1,2.828)(0,3)(1,2.828)(1.47,2.54976)%(1.45,2.5616)
   \psline[linewidth=0.33pt](1.53,2.51424)(2,2.236)(3,0)(-3,0)%(1.55,2.5024)(2,2.236)(3,0)(-3,0)
 \psline[linestyle=dotted,linewidth=0.33pt](0,0)(-2,2.236)
  \psline[linestyle=dotted,linewidth=0.33pt](0,0)(-1,2.828)
   \psline[linestyle=dotted,linewidth=0.33pt](0,0)(0,3)
    \psline[linestyle=dotted,linewidth=0.33pt](0,0)(1,2.828)
     \psline[linestyle=dotted,linewidth=0.33pt](0,0)(2,2.236)
        \rput(0,-0.2){\tiny $O$}
   }  }
%%%%%%%%%%%%%%%%%%%%%%%%%%%%%%%%%%%%%%%
   \rput{-114}(1.55,2.2){
  \scalebox{-0.15 -0.15}{
 \scalebox{0.7}{
 \rput(-7.5,0){
 \psline[linewidth=7pt](-1,0)(-8.55,3.775)
 \psline[linewidth=7pt](1,0)(14.35,6.675)
 \psecurve[linewidth=7pt](-8.5,3.5)(-8.55,3.775)(-9,9)(-8.5,15.5)(-7.6,20.33)(-7,23)
 \psecurve[linewidth=7pt](14.775,6.084)(14.35,6.675)(10.864,10.864)(6.084,14.775)(0,18)(-6.084,19.99)(-7.6,20.33)(-13.471,21.293)
        }}
        }}
        %%%%%%%%%%%%%%%%%%%%%%%%%%%%%%%%%%%%%%%
   \rput{-30.5}(10.1,5.3){
  \scalebox{-0.15 0.15}{
 \scalebox{0.7}{
 \rput(-7.5,0){
        \scalebox{0.7}{
\psline[linewidth=10pt](-1,0)(-8.55,3.775)
 \psline[linewidth=10pt](1,0)(14.35,6.675)
 \psecurve[linewidth=10pt](-8.5,3.5)(-8.55,3.775)(-9,9)(-8.5,15.5)(-7.6,20.33)(-7,23)
 \psecurve[linewidth=10pt](14.775,6.084)(14.35,6.675)(10.864,10.864)(6.084,14.775)(0,18)(-6.084,19.99)(-7.6,20.33)(-13.471,21.293)
                      }
        }}
        }}
 \rput(1.27,1.15){$x$}
 \rput(10.85,4.55){$x'$}
           \rput(-0.2,3.5){\Large $R_{ij}^{k_1}[x]$}
           \rput(12.5,6.6){\Large $R_{ji}^{k_2}[x']$}
 }

 \vspace{30mm}

Note that there exist values $\check{\Phi}_1 = \Phi_1(i,j)$ and
$\check{\Phi}_2 = \Phi_2(i,j)$ such that $-\pi/2 < \check{\Phi}_1 <
\check{\Phi}_2 < \pi/2$ and for any $x \in [P_{i-1},\ P_i]$,\
$\vphi_1(x)$ and $\vphi_2(x)$ belong to $[\check{\Phi}_1,\,
\check{\Phi}_2]$. If, in addition, $(i,j)$ does not coincide with
one of the pairs $(1,m)$,\, $(m,1)$ then there exists a value
$\check{\Phi}_0 = \Phi_0(i,j)$ such that $0 < \check{\Phi}_0 <
\check{\Phi}_2 - \check{\Phi}_1$ and for any $x \in [P_{i-1},\,
P_i]$,\, $|\vphi_2(x) - \vphi_1(x)| \ge \check{\Phi}_0$. Note in
passing that the measure preserving function $\vphi_{ij}$ was chosen
to be monotone decreasing just to make it possible to choose
$\check{\Phi}_1$ and $\check{\Phi}_2$; if it were defined as
monotone increasing then for $|i - j| = 1$ it would not be possible
to separate the functions $\vphi_1(x)$ and $\vphi_2(x)$ from
$-\pi/2$ and $\pi/2$. Note also that the condition $\s(1) \ne m$,
stated in lemma \ref{dense}, guarantees that the pair $(i,\s(i))$
does not coincide with one of the pairs $(1,m)$,\, $(m,1)$.

Denote by $\Om_{(i)} = \Om_{(i)}(r) = \{ \rho e_\theta:\, \rho \in
[0,\, r],\ \theta \in \Theta_i \}$ the circular sector bounded by
the radii $OP_{i-1}$ and $OP_i$. For any pair of values $i \ne j$,
the set $\Om_{(i,j)} = \Om_{(i,j)}(r,\del)$ is defined as union of
the triangle $OP_{i-1}P_i$ and a finite set of reflectors of the
form $R_{ij}^{k_l}[x_q^l]$, where $x_q^l \in [P_{i-1},\, P_i]$ and
$0 < k_l < 1$. The detailed definition is given in the two next
paragraphs. The reflectors do not mutually intersect and are
contained in the angle $\angle P_{i-1}OP_i := \{ \rho e_\theta:\,
\rho \ge 0,\ \theta \in \Theta_i \}$; besides, the set
$\pl\Om_{(i,j)} \cap [P_{i-1},\, P_i]$ is a union of segments of
common length at most 1. The last property can be formulated in a
slightly different manner: the collection of points of $[P_{i-1},\,
P_i]$ that do not belong to bases of the mentioned reflectors has
common length at most 1. The set $\Om = \Om(\s,r,\del)$ is defined
as the union
\begin{equation}\label{def Omega}
\Om = \left( \cup_{\s(i)=i} \Om_{(i)} \right) \cup \left(
\cup_{\s(i) \ne i} \Om_{(i,\s(i))} \right).
\end{equation}
Obviously, $\pl\Om$ contains $I$, and $\Om$ lies on one side of the
straight line $l_0$ containing $I$; hence $(\Om, I) \in \SSS$.

Let us give a detailed description of $\Om_{(i,j)}$. The trapezium
$\mathrm{T}(\del)$, defined in section \ref{opredelenie reflektora},
will be briefly referred to as $\mathrm{T}$. Denote by $\Lam_1^0$
the isosceles triangle with the base $[P_{i-1},\, P_i]$ and the
angle $\del$ at the base. Denote by $x_1^0$ the midpoint of the
segment $[P_{i-1},\, P_i]$ and consider the trapezium
$\mathrm{T}_1^0 := V_{x_1^0} k_0 \mathrm{T}$, where $k_0$ is the
greatest value, not exceeding 1, such that $\mathrm{T}_1^0$ is
contained in $\Lam_1^0$. The set $\Lam_1^0 \setminus \mathrm{T}_1^0$
is a union of three triangles; two of them have bases belonging to
$[P_{i-1},\, P_i]$; denote these triangles by $\Lam_1^1$ and
$\Lam_2^1$. Designate by $x_1^1$ and $x_2^1$ the midpoints of their
bases and consider trapezia $\mathrm{T}_1^1 = V_{x_1^1} k_1
\mathrm{T}$ and $\mathrm{T}_2^1 = V_{x_2^1} k_1 \mathrm{T}$, where
$k_1$ is the greatest value, not exceeding 1, such that
$\mathrm{T}_1^1 \subset \Lam_1^1$ and $\mathrm{T}_2^1 \subset
\Lam_2^1$.

Continuing this process, one obtains the sequence of numbers $k_l$,
triangles $\Lam_1^l, \ldots, \Lam_{2^l}^l$, points $x_1^l, \ldots,
x_{2^l}^l$, and trapezia $\mathrm{T}_1^l, \ldots,
\mathrm{T}_{2^l}^l$,\, $l = 0,\, 1,\, 2, \ldots$ (see the figure
below). These trapezia do not mutually intersect; the numbers $k_l$
form a decreasing geometric progression. The set $J_l = [P_{i-1},\,
P_i] \setminus \cup_{l'=0}^{l-1} \cup_{q=1}^{2^{l'}}
\mathrm{T}_q^{l'}$ is a union of $2^l$ intervals, besides the
one-dimensional Lebesgue measures $|J_l|$ of $J_l$,\, $l = 1,\,2,
\ldots$ form a decreasing geometric progression. Choose $l_0$ such
that $|J_{l_0}| < 1$. Denote by $\Lam_{(i)}$ the triangle
$OP_{i-1}P_i$ and let
$$
\Om_{(i,j)} = \Lam_{(i)} \cup \left( \cup_{l=0}^{l_0}
\cup_{q=1}^{2^l} R_{ij}^{k_l}[x_q^l] \right).
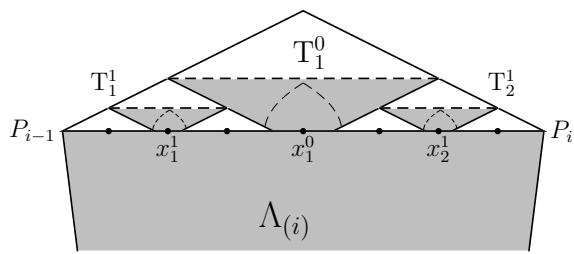
$$
 \vspace{-7mm}

\rput(6,-2.3){
       \scalebox{0.8}{
 \pspolygon(-4,0)(0,2)(4,0)
 \pspolygon[linecolor=white,fillstyle=solid,fillcolor=lightgray](-0.5,0)(-2.25,0.875)(2.25,0.875)(0.5,0)
       \psline(2.25,0.875)(0.5,0)(-0.5,0)(-2.25,0.875)
         \psline[linestyle=dashed](-2.25,0.875)(2.25,0.875)
           \pscurve[linewidth=0.35pt,linestyle=dashed](0,0.8)(0.15,0.7)(0.28,0.6)(0.39,0.5)(0.48,0.4)(0.6,0.2)(0.6336,0.08)
           \pscurve[linewidth=0.35pt,linestyle=dashed](0,0.8)(-0.15,0.7)(-0.28,0.6)(-0.39,0.5)(-0.48,0.4)(-0.6,0.2)(-0.6336,0.08)
 \pspolygon[linecolor=white,fillstyle=solid,fillcolor=lightgray](-3.75,-2)(-4,0)(4,0)(3.75,-2)
 \psline(-3.75,-2)(-4,0)(4,0)(3.75,-2)
  \rput(-4.5,-0){$P_{i-1}$}
  \rput(4.3,-0){$P_{i}$}
        \rput(0,-0.3){$x_1^0$}
        \rput(-2.25,-0.3){$x_1^1$}
        \rput(2.25,-0.3){$x_2^1$}
              \pspolygon[linecolor=white,fillstyle=solid,fillcolor=lightgray](2.03125,0)(1.265625,0.3828125)(3.234375,0.3828125)(2.46875,0)
              \psline(3.234375,0.3828125)(2.46875,0)(2.03125,0)(1.265625,0.3828125)
                                    \rput(2.25,0){
                                    \scalebox{0.4375}{
                                        \pscurve[linewidth=0.8pt,linestyle=dashed](0,0.8)(0.15,0.7)(0.28,0.6)(0.39,0.5)(0.48,0.4)(0.6,0.2)(0.6336,0.08)
                                        \pscurve[linewidth=0.8pt,linestyle=dashed](0,0.8)(-0.15,0.7)(-0.28,0.6)(-0.39,0.5)(-0.48,0.4)(-0.6,0.2)(-0.6336,0.08)
                                    }
                                    }
               \psline[linestyle=dashed](3.234375,0.3828125)(1.265625,0.3828125)
        \pspolygon[linecolor=white,fillstyle=solid,fillcolor=lightgray](-2.03125,0)(-1.265625,0.3828125)(-3.234375,0.3828125)(-2.46875,0)
        \psline(-3.234375,0.3828125)(-2.46875,0)(-2.03125,0)(-1.265625,0.3828125)
           \psline[linestyle=dashed](-3.234375,0.3828125)(-1.265625,0.3828125)
                                         \rput(-2.25,0){
                                         \scalebox{0.4375}{
                                        \pscurve[linewidth=0.8pt,linestyle=dashed](0,0.8)(0.15,0.7)(0.28,0.6)(0.39,0.5)(0.48,0.4)(0.6,0.2)(0.6336,0.08)
                                        \pscurve[linewidth=0.8pt,linestyle=dashed](0,0.8)(-0.15,0.7)(-0.28,0.6)(-0.39,0.5)(-0.48,0.4)(-0.6,0.2)(-0.6336,0.08)
                                         }
                                         }
   \psdots[dotsize=3pt](0,0)(2.25,0)(-2.25,0)(3.234375,0)(1.265625,0)(-3.234375,0)(-1.265625,0)
  \rput(0.1,1.3){\large $\mathrm{T}_{1}^0$}
     \rput(-0.3,-1.5){\Large $\Lam_{(i)}$}
   \rput(-3.3,0.8){$\mathrm{T}_{1}^1$}
   \rput(3.3,0.8){$\mathrm{T}_{2}^1$}
     }
  }
 \vspace{48mm}

At last, $\Om$ is defined by the formula (\ref{def Omega}). Each
reflector constituting the subset $\Om_{(i,\s(i))}$ will be called
{\it $i$-reflector} of the set $\Om$.

On the figure below, the pair of sets $(\Om = \Om(\s,r,\del),\ I)$
is shown, where $\s = \big( {\scriptsize
\begin{array}{ccccc} 1 & 2 & 3 & 4 & 5\\ 1 & 4 & 3 & 2 & 5
\end{array} \big)}$,\, $r = 2.5.$,\, $\del = 0.4$.

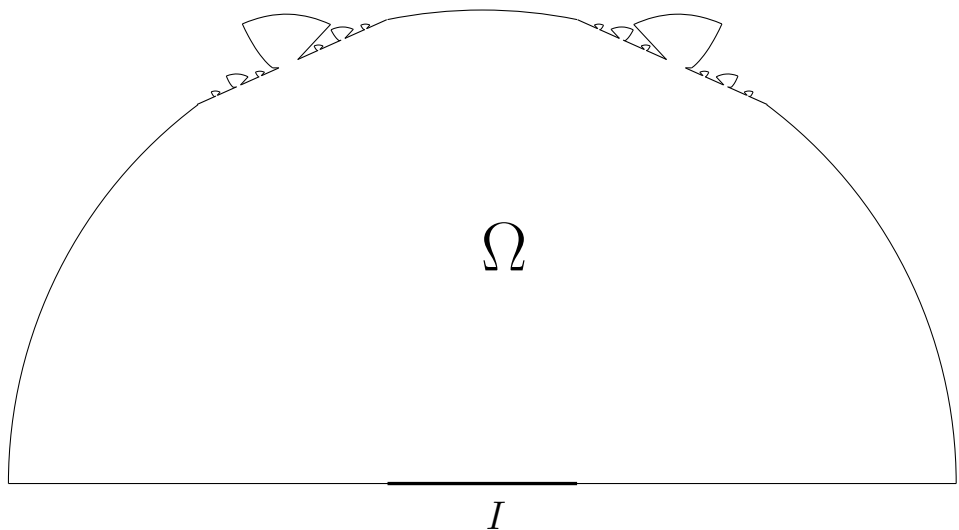
\vspace*{75mm}

  \rput(6.5,0){
  \scalebox{1.05}{
  %\pscircle(0,0){6}
  %\psdots(-6,0)(-3.6,4.8)(-1.2,5.878775)(1.2,5.878775)(3.6,4.8)(6,0)
  %\psline[linewidth=0.38pt](-1.2,5.878775)(1.2,5.878775)
  %\psline[linewidth=0.38pt](-3.6,4.8)(-6,0)(6,0)(3.6,4.8)
      \psline[linewidth=0.38pt](-6,0)(6,0)
      \psline[linewidth=1.2pt](-1.2,0)(1.2,0)
  \psarc[linewidth=0.38pt](0,0){6}{0}{53.13}
  \psarc[linewidth=0.38pt](0,0){6}{78.463}{101.537}
  \psarc[linewidth=0.38pt](0,0){6}{126.87}{180}
 %%%%%%%%%%%%%%%%%%%%%%%%%%%%%%%%%%%%%%%%%%%% right part
 %%%%%%%%%%%%%%%%%%%%%%%%%%%%%%%%%%%%%%%%%%%% right part
    \rput{-24.1}(2.43,5.33){
    \scalebox{1.316}{
         \psline[linewidth=0.3pt](-1,0)(-0.793,0)
          \psline[linewidth=0.3pt](-0.77,0)(-0.56,0)
        \psline[linewidth=0.3pt](-0.52,0)(-0.30,0)
         \psline[linewidth=0.3pt](-0.277,0)(-0.06,0)
        \psline[linewidth=0.3pt](0.14,0)(0.318,0)
        \psline[linewidth=0.3pt](0.345,0)(0.535,0)
          \psline[linewidth=0.3pt](0.58,0)(0.78,0)
          \psline[linewidth=0.3pt](0.80,0)(1,0)
         %\psdot[dotsize=0.5pt](0.55,0)
                  %\psline[linestyle=dashed,linewidth=0.1pt](-1,0)(-0.08,0.39)
                  %\psline[linestyle=dashed,linewidth=0.1pt](1,0)(0.08,0.39)
                                 % 11111111111111111----------------------------
                             \rput(-0.78,0){
        \scalebox{0.0256}{
  %\psline[linewidth=0.1pt](0,0)(0.734,0.321)
  %%%%                          broken line:
  %\psline[linewidth=4pt](0.4,0)(-0.4,0)
  \psline[linewidth=11pt](0.4,0)(0.625,0.098)
  \psline[linewidth=11pt](-0.4,0)(-2.17,0.7503)
               %\psline[linestyle=dashed,linewidth=0.1pt](-0.4,0)(-2.24,0.78)
      %\psline[linestyle=dashed,linewidth=0.1pt](0.4,0)(2.24,0.78)
%%%%%%%%%%%%%%%%%%%%%%%%%%%%%%%% right parabola
 %\psdot[dotsize=20pt](0,0)
          \rput{-2}(-0.14,0){
  \psecurve[linewidth=11pt](0.38,-0.283)(0.625,0.098)(0.732,0.32)(0.89,0.89)(0.979,1.6)(0.999,2.12)(0.999,2.6) % right curve
%%%%%%%%%%%%%%%%%%%%%%%%%%%%%%%%%
 \psecurve[linewidth=11pt](-2.25,0.5)(-2.117,0.681)(-1.5,1.28)(-0.8,1.725)(0,2.005)(0.999,2.12)(2,2.04) % left curve
          }
          }
          }
                     % 1111111111111111111111111111111111111111111111111111111111111111111111
                             \rput(-0.54,0){
        \scalebox{0.0625}{
  %\psline[linewidth=0.1pt](0,0)(0.734,0.321)
  %%%%                          broken line:
  %\psline[linewidth=1.6pt](0.4,0)(-0.4,0)
  \psline[linewidth=4.6pt](0.4,0)(0.625,0.098)
  \psline[linewidth=4.6pt](-0.4,0)(-2.17,0.7503)
               %\psline[linestyle=dashed,linewidth=0.1pt](-0.4,0)(-2.24,0.78)
      %\psline[linestyle=dashed,linewidth=0.1pt](0.4,0)(2.24,0.78)
%%%%%%%%%%%%%%%%%%%%%%%%%%%%%%%% right parabola
 %\psdot[dotsize=6.25pt](0,0)
          \rput{-2}(-0.14,0){
  \psecurve[linewidth=4.6pt](0.38,-0.283)(0.625,0.098)(0.732,0.32)(0.89,0.89)(0.979,1.6)(0.999,2.12)(0.999,2.6) % right curve
%%%%%%%%%%%%%%%%%%%%%%%%%%%%%%%%%
 \psecurve[linewidth=4.6pt](-2.25,0.5)(-2.117,0.681)(-1.5,1.28)(-0.8,1.725)(0,2.005)(0.999,2.12)(2,2.04) % left curve
          }
          }
          }
        % 22222222222222222222222222222222----------------------------------
        \rput(-0.285,0){
        \scalebox{0.0256}{
  %\psline[linewidth=0.4pt](0,0)(0.734,0.321)
  %%%%                          broken line:
  %\psline[linewidth=4pt](0.4,0)(-0.4,0)
  \psline[linewidth=11pt](0.4,0)(0.625,0.098)
  \psline[linewidth=11pt](-0.4,0)(-2.09,0.715)
               %\psline[linestyle=dashed,linewidth=0.4pt](-0.4,0)(-2.24,0.78)
                %\psline[linestyle=dashed,linewidth=0.4pt](0.4,0)(2.24,0.78)
 %%%%%%%%%%%%%%%%%%%%%%%%%%%%%%%%
 %\psdot[dotsize=20pt](0,0)
  \psecurve[linewidth=11pt](0.38,-0.283)(0.625,0.098)(0.732,0.32)(0.89,0.89)(0.979,1.6)(0.999,2.12)(0.999,2.6) % right curve
 \psecurve[linewidth=11pt](-2.25,0.5)(-2.09,0.715)(-1.5,1.28)(-0.8,1.725)(0,2.005)(0.999,2.12)(2,2.04) % left curve
             }
              }
      % 2222222222222222222222222222222222222222222222222222222222222222222222222
        \scalebox{0.25}{
  %\psline[linewidth=0.4pt](0,0)(0.734,0.321)
  %%%%                          broken line:
  %\psline[linewidth=0.4pt](0.4,0)(-0.4,0)
  \psline[linewidth=1.16pt](0.4,0)(0.625,0.098)
  \psline[linewidth=1.16pt](-0.4,0)(-2.09,0.715)
               %\psline[linestyle=dashed,linewidth=0.4pt](-0.4,0)(-2.24,0.78)
                %\psline[linestyle=dashed,linewidth=0.4pt](0.4,0)(2.24,0.78)
 %%%%%%%%%%%%%%%%%%%%%%%%%%%%%%%%
 %\psdot[dotsize=2pt](0,0)
  \psecurve[linewidth=1.16pt](0.38,-0.283)(0.625,0.098)(0.732,0.32)(0.89,0.89)(0.979,1.6)(0.999,2.12)(0.999,2.6) % right curve
 \psecurve[linewidth=1.16pt](-2.25,0.5)(-2.09,0.715)(-1.5,1.28)(-0.8,1.725)(0,2.005)(0.999,2.12)(2,2.04) % left curve
             }
       % 22222222222222222222222222222222+++++++++++++++++++++++++++++++
        \rput(0.16,0){
        \scalebox{0.0256}{
  %\psline[linewidth=0.4pt](0,0)(0.734,0.321)
  %%%%                          broken line:
  %\psline[linewidth=4pt](0.4,0)(-0.4,0)
  \psline[linewidth=11pt](0.4,0)(0.625,0.098)
  \psline[linewidth=11pt](-0.4,0)(-2.09,0.715)
               %\psline[linestyle=dashed,linewidth=0.4pt](-0.4,0)(-2.24,0.78)
                %\psline[linestyle=dashed,linewidth=0.4pt](0.4,0)(2.24,0.78)
 %%%%%%%%%%%%%%%%%%%%%%%%%%%%%%%%
 %\psdot[dotsize=20pt](0,0)
  \psecurve[linewidth=11pt](0.38,-0.283)(0.625,0.098)(0.732,0.32)(0.89,0.89)(0.979,1.6)(0.999,2.12)(0.999,2.6) % right curve
 \psecurve[linewidth=11pt](-2.25,0.5)(-2.09,0.715)(-1.5,1.28)(-0.8,1.725)(0,2.005)(0.999,2.12)(2,2.04) % left curve
             }
              }
                % 333333333333333333333333333333333333333333333333333333333333333333333
        \rput(0.39,0){
        \scalebox{0.0625}{
%%%%%%%%%       broken line
  %\psline[linewidth=1.6pt](0.4,0)(-0.4,0)
  \psline[linewidth=4.6pt](0.4,0)(0.76,0.15)
         %\psline[linewidth=1.6pt,linestyle=dashed](0.4,0)(1.32,0.39)
  \psline[linewidth=4.6pt](-0.4,0)(-2.0,0.68)
       % \psline[linewidth=1.6pt,linestyle=dashed](-0.4,0)(-2.1,0.7206)
%%%%%%%%%%%%%%%%%%%%%%%%%%%%%%%%%%%%%%%%%%%%%%%%
 %\psdots[dotsize=6.25pt](0,0)%(0,2)(0.888,0.54)
 %%%%%%%%%%%%%%%%%%%%%%%%%%%%%%%%
 \psecurve[linewidth=4.6pt](0.57,-0.2)(0.76,0.15)(0.89,0.54)(0.985,1.09)(1.01,1.6)(1.005,2.16)(0.99,2.42)  % right curve
 \psecurve[linewidth=4.6pt](-2.33,0.2)(-2.0,0.68)(-1.46,1.235)(-0.78,1.7)(0,2.005)(1.005,2.16)(2,2.125) %dots of left curve
  }
  }
                  % 33333333333333333333333333333++++++++++++++++++++
        \rput(0.625,0){
        \scalebox{0.0256}{
%%%%%%%%%       broken line
  %\psline[linewidth=4pt](0.4,0)(-0.4,0)
  \psline[linewidth=11pt](0.4,0)(0.76,0.15)
         %\psline[linewidth=4pt,linestyle=dashed](0.4,0)(1.32,0.39)
  \psline[linewidth=11pt](-0.4,0)(-2.0,0.68)
       % \psline[linewidth=1.6pt,linestyle=dashed](-0.4,0)(-2.1,0.7206)
%%%%%%%%%%%%%%%%%%%%%%%%%%%%%%%%%%%%%%%%%%%%%%%%
 %\psdots[dotsize=20pt](0,0)%(0,2)(0.888,0.54)
 %%%%%%%%%%%%%%%%%%%%%%%%%%%%%%%%
 \psecurve[linewidth=11pt](0.57,-0.2)(0.76,0.15)(0.89,0.54)(0.985,1.09)(1.01,1.6)(1.005,2.16)(0.99,2.42)  % right curve
 \psecurve[linewidth=11pt](-2.33,0.2)(-2.0,0.68)(-1.46,1.235)(-0.78,1.7)(0,2.005)(1.005,2.16)(2,2.125) %dots of left curve
  }
  }
    }
    }
%%%%%%%%%%%%%%%%%%%%%%%%%%%%%%%%%%%%%%%%%%%% left part
%%%%%%%%%%%%%%%%%%%%%%%%%%%%%%%%%%%%%%%%%%%% left part
    \scalebox{-1 1}{
    \rput{-24.1}(2.43,5.33){
    \scalebox{1.316}{
         \psline[linewidth=0.3pt](-1,0)(-0.793,0)
          \psline[linewidth=0.3pt](-0.77,0)(-0.56,0)
        \psline[linewidth=0.3pt](-0.52,0)(-0.30,0)
         \psline[linewidth=0.3pt](-0.277,0)(-0.06,0)
        \psline[linewidth=0.3pt](0.14,0)(0.318,0)
        \psline[linewidth=0.3pt](0.345,0)(0.535,0)
          \psline[linewidth=0.3pt](0.58,0)(0.78,0)
          \psline[linewidth=0.3pt](0.80,0)(1,0)
         %\psdot[dotsize=0.5pt](0.55,0)
                  %\psline[linestyle=dashed,linewidth=0.1pt](-1,0)(-0.08,0.39)
                  %\psline[linestyle=dashed,linewidth=0.1pt](1,0)(0.08,0.39)
                                 % 11111111111111111----------------------------
                             \rput(-0.78,0){
        \scalebox{0.0256}{
  %\psline[linewidth=11pt](0,0)(0.734,0.321)
  %%%%                          broken line:
  %\psline[linewidth=11pt](0.4,0)(-0.4,0)
  \psline[linewidth=11pt](0.4,0)(0.625,0.098)
  \psline[linewidth=11pt](-0.4,0)(-2.17,0.7503)
               %\psline[linestyle=dashed,linewidth=0.1pt](-0.4,0)(-2.24,0.78)
      %\psline[linestyle=dashed,linewidth=0.1pt](0.4,0)(2.24,0.78)
%%%%%%%%%%%%%%%%%%%%%%%%%%%%%%%% right parabola
 %\psdot[dotsize=20pt](0,0)
          \rput{-2}(-0.14,0){
  \psecurve[linewidth=11pt](0.38,-0.283)(0.625,0.098)(0.732,0.32)(0.89,0.89)(0.979,1.6)(0.999,2.12)(0.999,2.6) % right curve
%%%%%%%%%%%%%%%%%%%%%%%%%%%%%%%%%
 \psecurve[linewidth=11pt](-2.25,0.5)(-2.117,0.681)(-1.5,1.28)(-0.8,1.725)(0,2.005)(0.999,2.12)(2,2.04) % left curve
          }
          }
          }
                     % 1111111111111111111111111111111111111111111111111111111111111111111111
                             \rput(-0.54,0){
        \scalebox{0.0625}{
  %\psline[linewidth=0.1pt](0,0)(0.734,0.321)
  %%%%                          broken line:
  %\psline[linewidth=1.6pt](0.4,0)(-0.4,0)
  \psline[linewidth=4.6pt](0.4,0)(0.625,0.098)
  \psline[linewidth=4.6pt](-0.4,0)(-2.17,0.7503)
               %\psline[linestyle=dashed,linewidth=0.1pt](-0.4,0)(-2.24,0.78)
      %\psline[linestyle=dashed,linewidth=0.1pt](0.4,0)(2.24,0.78)
%%%%%%%%%%%%%%%%%%%%%%%%%%%%%%%% right parabola
 %\psdot[dotsize=6.25pt](0,0)
          \rput{-2}(-0.14,0){
  \psecurve[linewidth=4.6pt](0.38,-0.283)(0.625,0.098)(0.732,0.32)(0.89,0.89)(0.979,1.6)(0.999,2.12)(0.999,2.6) % right curve
%%%%%%%%%%%%%%%%%%%%%%%%%%%%%%%%%
 \psecurve[linewidth=4.6pt](-2.25,0.5)(-2.117,0.681)(-1.5,1.28)(-0.8,1.725)(0,2.005)(0.999,2.12)(2,2.04) % left curve
          }
          }
          }
        % 22222222222222222222222222222222----------------------------------
        \rput(-0.285,0){
        \scalebox{0.0256}{
  %\psline[linewidth=0.4pt](0,0)(0.734,0.321)
  %%%%                          broken line:
  %\psline[linewidth=4pt](0.4,0)(-0.4,0)
  \psline[linewidth=11pt](0.4,0)(0.625,0.098)
  \psline[linewidth=11pt](-0.4,0)(-2.09,0.715)
               %\psline[linestyle=dashed,linewidth=0.4pt](-0.4,0)(-2.24,0.78)
                %\psline[linestyle=dashed,linewidth=0.4pt](0.4,0)(2.24,0.78)
 %%%%%%%%%%%%%%%%%%%%%%%%%%%%%%%%
 %\psdot[dotsize=20pt](0,0)
  \psecurve[linewidth=11pt](0.38,-0.283)(0.625,0.098)(0.732,0.32)(0.89,0.89)(0.979,1.6)(0.999,2.12)(0.999,2.6) % right curve
 \psecurve[linewidth=11pt](-2.25,0.5)(-2.09,0.715)(-1.5,1.28)(-0.8,1.725)(0,2.005)(0.999,2.12)(2,2.04) % left curve
             }
              }
      % 2222222222222222222222222222222222222222222222222222222222222222222222222
        \scalebox{0.25}{
  %\psline[linewidth=0.4pt](0,0)(0.734,0.321)
  %%%%                          broken line:
  %\psline[linewidth=0.4pt](0.4,0)(-0.4,0)
  \psline[linewidth=1.16pt](0.4,0)(0.625,0.098)
  \psline[linewidth=1.16pt](-0.4,0)(-2.09,0.715)
               %\psline[linestyle=dashed,linewidth=0.4pt](-0.4,0)(-2.24,0.78)
                %\psline[linestyle=dashed,linewidth=0.4pt](0.4,0)(2.24,0.78)
 %%%%%%%%%%%%%%%%%%%%%%%%%%%%%%%%
 %\psdot[dotsize=2pt](0,0)
  \psecurve[linewidth=1.16pt](0.38,-0.283)(0.625,0.098)(0.732,0.32)(0.89,0.89)(0.979,1.6)(0.999,2.12)(0.999,2.6) % right curve
 \psecurve[linewidth=1.16pt](-2.25,0.5)(-2.09,0.715)(-1.5,1.28)(-0.8,1.725)(0,2.005)(0.999,2.12)(2,2.04) % left curve
             }
       % 22222222222222222222222222222222+++++++++++++++++++++++++++++++
        \rput(0.16,0){
        \scalebox{0.0256}{
  %\psline[linewidth=0.4pt](0,0)(0.734,0.321)
  %%%%                          broken line:
  %\psline[linewidth=4pt](0.4,0)(-0.4,0)
  \psline[linewidth=11pt](0.4,0)(0.625,0.098)
  \psline[linewidth=11pt](-0.4,0)(-2.09,0.715)
               %\psline[linestyle=dashed,linewidth=0.4pt](-0.4,0)(-2.24,0.78)
                %\psline[linestyle=dashed,linewidth=0.4pt](0.4,0)(2.24,0.78)
 %%%%%%%%%%%%%%%%%%%%%%%%%%%%%%%%
 %\psdot[dotsize=20pt](0,0)
  \psecurve[linewidth=11pt](0.38,-0.283)(0.625,0.098)(0.732,0.32)(0.89,0.89)(0.979,1.6)(0.999,2.12)(0.999,2.6) % right curve
 \psecurve[linewidth=11pt](-2.25,0.5)(-2.09,0.715)(-1.5,1.28)(-0.8,1.725)(0,2.005)(0.999,2.12)(2,2.04) % left curve
             }
              }
                % 333333333333333333333333333333333333333333333333333333333333333333333
        \rput(0.39,0){
        \scalebox{0.0625}{
%%%%%%%%%       broken line
  %\psline[linewidth=1.6pt](0.4,0)(-0.4,0)
  \psline[linewidth=4.6pt](0.4,0)(0.76,0.15)
         %\psline[linewidth=1.6pt,linestyle=dashed](0.4,0)(1.32,0.39)
  \psline[linewidth=4.6pt](-0.4,0)(-2.0,0.68)
       % \psline[linewidth=1.6pt,linestyle=dashed](-0.4,0)(-2.1,0.7206)
%%%%%%%%%%%%%%%%%%%%%%%%%%%%%%%%%%%%%%%%%%%%%%%%
 %\psdots[dotsize=6.25pt](0,0)%(0,2)(0.888,0.54)
 \rput{211}(-0.73,1.73){
 \scalebox{3.7}{
 }
 }
 %%%%%%%%%%%%%%%%%%%%%%%%%%%%%%%%
 \psecurve[linewidth=4.6pt](0.57,-0.2)(0.76,0.15)(0.89,0.54)(0.985,1.09)(1.01,1.6)(1.005,2.16)(0.99,2.42)  % right curve
 \psecurve[linewidth=4.6pt](-2.33,0.2)(-2.0,0.68)(-1.46,1.235)(-0.78,1.7)(0,2.005)(1.005,2.16)(2,2.125) %dots of left curve
  }
  }
                  % 33333333333333333333333333333++++++++++++++++++++
        \rput(0.625,0){
        \scalebox{0.0256}{
%%%%%%%%%       broken line
  %\psline[linewidth=4pt](0.4,0)(-0.4,0)
  \psline[linewidth=11pt](0.4,0)(0.76,0.15)
         %\psline[linewidth=4pt,linestyle=dashed](0.4,0)(1.32,0.39)
  \psline[linewidth=11pt](-0.4,0)(-2.0,0.68)
       % \psline[linewidth=1.6pt,linestyle=dashed](-0.4,0)(-2.1,0.7206)
%%%%%%%%%%%%%%%%%%%%%%%%%%%%%%%%%%%%%%%%%%%%%%%%
 %\psdots[dotsize=20pt](0,0)%(0,2)(0.888,0.54)
 \rput{211}(-0.73,1.73){
 \scalebox{3.7}{
 }
 }
 %%%%%%%%%%%%%%%%%%%%%%%%%%%%%%%%
 \psecurve[linewidth=11pt](0.57,-0.2)(0.76,0.15)(0.89,0.54)(0.985,1.09)(1.01,1.6)(1.005,2.16)(0.99,2.42)  % right curve
 \psecurve[linewidth=11pt](-2.33,0.2)(-2.0,0.68)(-1.46,1.235)(-0.78,1.7)(0,2.005)(1.005,2.16)(2,2.125) %dots of left curve
  }
  }
    }
    }
    }
 \rput(-0.1,-0.4){\large $I$}
 \rput(0.0,3){\Huge $\Om$}

  }
  }

\vspace{15mm}

In the particular case where $\s =$ id, that is, $\s(i) = i$ for any
$i$, the set $\Om = \Om(\text{id},r,\del)$ takes an especially
simple form: it is just the upper half-circle of radius $r$.

Note that there exist values $\Phi_0$,\, $\Phi_1$,\, $\Phi_2$, not
depending on $r$ and $\del$, such that $-\pi/2 < \Phi_1 < \Phi_2 <
\pi/2$,\ $0 < \Phi_0 < \Phi_2 - \Phi_1$, and for any
$(\vphi_1,\vphi_2,\del)$-reflector, being a constituent of $\Om$,
holds $\vphi_1,\, \vphi_2 \in [\Phi_1,\, \Phi_2]$,\ $|\vphi_2 -
\vphi_1| \ge \Phi_0$. It suffices to put $\Phi_1 = \min_{i \ne j}
\Phi_1(i,j)$, $\Phi_2 = \max_{i \ne j} \Phi_2(i,j)$,\ $\Phi_0 = \min
\{ \Phi_0(i,j):\, i \ne j,\ (i,j) \ne (1,m),\, (m,1) \}$.

From existence of the mentioned values and from lemma \ref{reflector
lemma} it follows

\begin{corollary}\label{sledstvie c}
If $x$,\, $\vphi$ are such that $x$ belongs to the $(1 - c\del
-c|\vphi - \vphi_1|)$-base of some
$(\vphi_1,\vphi_2,\del)$-reflector constituting $\Om(\s,r,\del)$
then
\begin{equation}\label{sl phi}
\vphi^+_{\check{R},\check{I}} (x, \vphi) - \vphi_2 =
-\frac{\cos\vphi_1}{\cos\vphi_2}\, (\vphi - \vphi_1) + O(\del) +
O((\vphi - \vphi_1)^2).
\end{equation}
Here $\check{R}$ denotes this reflector, $\check{I}$, its base. If,
in addition, $\del < c_0$ and $|\vphi - \vphi_1| < c_0$ then
\begin{equation}\label{sl n}
n_{\check{R},\check{I}} (x, \vphi) = 3.
\end{equation}
Here $c$,\, $c_0$ are positive values that depend only on $\s$.
\end{corollary}

Like in lemma \ref{reflector lemma}, here we denote by $O(f)$ a
generic function of $x,\, \vphi$, $\vphi_1,\, \vphi_2$, $\del,\, r$
such that $O(f)/f$ is bounded over all admissible values of the
arguments.

\subsection{Proof of convergence $\nu_{\Om,I} \to \nu^\s$}

Consider a positive function $\del = \del(r) = o(1/r)$,\, $r \to
+\infty$. Recall that $n_{\Om,I}(x,\vphi)$ designates the number of
reflections from $\pl\Om$ (including the last reflection from $I$)
the particle getting out of $x \in I$ in the direction $\vphi$ has
to make before returning to $I$. Define the function $n_\s:\, I
\times \interval \to \RRR$ by the relations~ $n_\s(x,\vphi) = 5$, if
$\vphi \in \Theta_i$,\, $\s(i) \ne i$, and $n_\s(x,\vphi) = 2$, if
$\vphi \in \Theta_i$,\, $\s(i) = i$. Thus, $n_\s$ does not depend on
the argument $x$:\ $n_\s(x,\vphi) = n_\s(\vphi)$. At the endpoints
of $\Theta_i$, the definition may turn out to be ambiguous; correct
it by redefining $n_\s$ at these points. Denote $\Om(r) :=
\Om(\s,r,\del(r))$; that is, $\Om(r)$ is a one parameter family of
sets depending on $r$. Recall that $(\Om(r), I) \in \SSS$.

The following lemma is basic.

\begin{lem}\label{l4}
One has convergence in measure
\begin{equation}\label{conv n}
\lim_{r\to+\infty} n_{\Om(r),I} = n_\s,
\end{equation}
\begin{equation}\label{conv fi}
\lim_{r\to+\infty} \vphi^+_{\Om(r),I} = \vphi_\s.
\end{equation}
\end{lem}

In other words, the lemma states that for large $r$, all the
particles getting in $\Om$ through $I$, except for a small part of
them, after one or four reflections get out through $I$, besides the
relation between the angles of getting in and getting out is close
to the one given by $\vphi_\s$.

\begin{proof}
Fix $i$. It suffices to prove that there exists a family of sets
$A_r \subset I \times \Theta_i$ such that their measures converge to
zero: $\lim_{r \to +\infty} \mu(A_r) = 0$ and
\begin{equation}\label{convergence n}
\text{for} \ \ (x,\vphi) \in (I \times \Theta_i) \setminus A_r \ \
\text{holds} \ \ n_{\Om(r),I}(x,\vphi) = n_\s(\vphi)
\end{equation}
and
\begin{equation}\label{convergence fi}
\lim_{r\to +\infty} \sup_{(x,\vphi) \in (I \times \Theta_i)
\setminus A_r} |\vphi^+_{\Om(r),I}(x,\vphi) - \vphi_\s(\vphi)| = 0.
\ \ \ \ \ \ \ \ \
\end{equation}
Consider the cases $\s(i) = i$ and $\s(i) \ne i$ separately.
 \vspace{2mm}

(a) Let $\s(i) = i$. In this case the relations (\ref{convergence
n}) and (\ref{convergence fi}) may be re-written in the form
\begin{equation*}\label{convergence n: s(i)=i}
\text{1)} \ \ \text{for} \ \ (x,\vphi) \in (I \times \Theta_i)
\setminus A_r \ \ \text{holds} \ \ n_{\Om(r),I}(x,\vphi) = 2; \ \ \
\ \ \ \ \ \ \ \ \ \ \ \ \ \
\end{equation*}
\begin{equation*}\label{convergence fi: s(i)not=i}
\text{2)} \ \ \lim_{r\to +\infty} \sup_{(x,\vphi) \in (I \times
\Theta_i) \setminus A_r} |\vphi^+_{\Om(r),I}(x,\vphi) - \vphi| = 0.
\ \ \ \ \ \ \ \ \ \ \ \ \ \ \ \ \ \ \ \ \ \ \ \ \ \ \ \
\end{equation*}
It follows from simple geometric argument that if a billiard
particle in $\Om(r)$ is initially placed at $x = (\xi, 0) \in I$ and
has velocity $e_\vphi$,\, $\vphi \in (\theta_{i-1} +
\arcsin\frac{1}{2r},\, \theta_i - \arcsin\frac{1}{2r}) \subset
\Theta_i$ then the first reflection point $x' = r e_\theta$ belongs
to the arc $\stackrel{\smile}{P_{i-1} P_i} = \{ r e_\theta: \theta
\in \Theta_i \}$. The second reflection occurs at a point $x^+ =
(\xi^+, 0)$, the velocity before the second reflection being
$-e_{\vphi^+}$. By using simple geometric argument, related to the
triangle with vertices $x$,\, $x'$, and $x^+$, one comes to the
relation
$$
\frac{1}{\xi} - \frac{1}{\xi^+} = \frac{2\sin\theta}{r}\,.
$$
It implies that if $|\xi| < 1/(2 + 4/r)$ then $|\xi^+| < 1/2$, that
is, the second reflection point $x^+$ belongs to $I$. Denoting
$$
A_r = \big( \big[-\frac{1}{2},\, -\frac{1}{2 + {4}/{r}}\big] \cup
\big[\frac{1}{2 + {4}/{r}},\, \frac{1}{2}\big] \big) \times \big(
\big[\theta_{i-1},\, \theta_{i-1} + \arcsin\frac{1}{2r}\big] \cup
\big[\theta_i - \arcsin\frac{1}{2r},\, \theta_i\big] \big),
$$
one concludes that for $(x,\vphi) \in (I \times \Theta_i) \setminus
A_r$,~ 1) holds true; besides, one has $\vphi^+ =
\vphi^+_{\Om(r),I}(x,\vphi)$. One obviously has $\lim_{r \to
+\infty} \mu(A_r) = 0$.

 
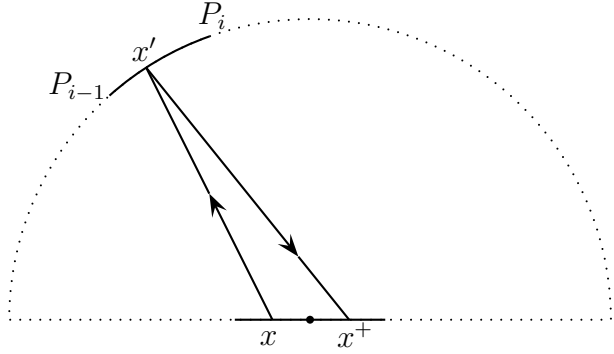
\vspace*{50mm}

 \psarc[linestyle=dotted,linewidth=0.8pt](6,0){4}{0}{180}
 \psarc[linewidth=0.8pt](6,0){4}{109.5}{131.8}
 \psline[linestyle=dotted,linewidth=0.8pt](2,0)(10,0)
 \psline[linewidth=0.8pt](5,0)(7,0)
     \psline[linewidth=0.8pt,arrows=->,arrowscale=2](5.5,0)(4.6607,1.67735)
     \psline[linewidth=0.8pt](4.6607,1.67735)(3.8214,3.3547)
  \psline[linewidth=0.8pt](6.52,0)(5.84535,0.838675)
     \psline[linewidth=0.8pt,arrows=<-,arrowscale=2](5.84535,0.838675)(3.8214,3.3547)
    \psdot[dotsize=3pt](6,0)
      \rput(2.9,3.1){$P_{i-1}$}
      \rput(4.7,4.05){$P_{i}$}
      \rput(3.8,3.6){$x'$}
       \rput(5.45,-0.25){$x$}
      \rput(6.6,-0.2){$x^+$}

 \vspace{10mm}

Further, from geometric argument related to the same triangle it
follows that $|\vphi^+ - \vphi| < 2\arctan \frac{1}{2r} = o(1)$,\,
$r \to +\infty$. Thus, 2) is also proved.
 \vspace{2mm}

(b) Let $\s(i) = j \ne i$. All the estimates below are true for $r
\to +\infty$ and are uniform with respect to $x$ and $\vphi$. The
assertions stated here are true for sufficiently large values of
$r$, and therefore for sufficiently small $\del = \del(r)$.

Let $\phi(r) := \sup \measuredangle LPO$, the supremum being taken
over all points $L \in I$ and over all points $P$ contained in the
broken line $P_0 P_1 \ldots P_m$. Denote $\kappa(r) = 1 - c\del(r) -
c\phi(r)$, where $c$ is the constant defined in corollary
\ref{sledstvie c}. Obviously, $\kappa(r) = 1 + O(r)$. Denote by
$[P_{i-1},\, P_i]'$ the union of $\kappa(r)$-bases of the
$i$-reflectors that have centers located at distances more than 1
from both points $P_{i-1}$,\, $P_i$. One has $[P_{i-1},\, P_i]'
\subset [P_{i-1},\, P_i]$, besides the length of $[P_{i-1},\, P_i]
\setminus [P_{i-1},\, P_i]'$ is $O(1)$. The set $[P_{j-1},\, P_j]'$
is defined analogously.

Consider the set $A^{(1)}_r$ of values $(x,\vphi) \in I \times
\Theta_i$ such that the ray, getting out of the point $x$ in the
direction $\vphi$, intersects the broken line $P_0 P_1 \ldots P_m$
at a point that does not belong to $[P_{i-1},\, P_i]'$. Thus, the
point of intersection $x'$ belongs to a set of Lebesgue measure
$O(1)$. Denote by $\vphi'$ the angle the vector $e_\vphi$ forms with
$n_{i}$. (Recall that $n_{i}$ is the unit outer normal to $\pl\QQQ$
at a point of $[P_{i-1},\, P_i]$). For any fixed $x'$, the value
$\vphi'$ is contained in an interval of length $O(1/r)$. Thus,
measure $\mu$ of the set of points $(x',\vphi')$ is $O(1/r)$. Taking
into account that the mapping $(x,\vphi) \to (x',\vphi')$ preserves
the measure $\mu$, one concludes that $\mu(A^{(1)}_r) = O(1/r)$,\,
$r \to +\infty$.

A billiard particle with initial data $(x,\vphi) \in (I \times
\Theta_i) \setminus A^{(1)}_r$ at some moment intersects the
$\kappa(r)$-base of an $i$-reflector, makes two reflections off its
boundary, then, according to corollary \ref{sledstvie c}, intersects
the base again and leaves the reflector. At the moment of leaving,
its velocity forms an angle $O(1/r)$ with the ray $AB$ that is
defined as follows: the point $A$ is the center of this
$i$-reflector, $\theta$ is such that $\frac{OA}{|OA|} = e_\theta$;\
$\theta^+ = \vphi_\s(\theta)$, the point $B$ belongs to the broken
line $P_0 P_1 \ldots P_m$ and $\frac{OB}{|OB|} = e_{\theta^+}$.

Let $x''$ be the point at which the particle intersects the line
$\pl\QQQ = P_0 P_1 \ldots P_m P_0$ again (for the third time), and
let $\vphi''$ be the angle the velocity of the particle forms with
$n_j$ at the moment of intersection. Denote by $A^{(2)}_r$ the set
of values $(x,\vphi) \in (I \times \Theta_i) \setminus A^{(1)}_r$
such that $x'' \not\in [P_{j-1},\, P_j]'$. Again, it is true that
$x''$ belongs to a set of Lebesgue measure $O(1)$ and for any fixed
$x''$, the value $\vphi''$ belongs to an interval of length
$O(1/r)$. It follows that measure of the set of points
$(x'',\vphi'')$ is $O(1/r)$, and taking into account that the
mapping $(x,\vphi) \to (x'',\vphi'')$ preserves the measure $\mu$,
one concludes that $\mu(A^{(2)}_r) = O(1/r)$,\, $r \to +\infty$.

Consider an element $(x,\vphi) \in (I \times \Theta_i) \setminus
(A^{(1)}_r \cup A^{(2)}_r)$; the corresponding billiard particle
reflects two times off the boundary of an $i$-reflector, then
reflects two times off the boundary of a $j$-reflector, and for the
fifth time, reflects at a point $x^+ \in [-r,\, r] \times \{ 0 \}$.
Denote by $-e_{\vphi^+}$ the velocity just before the fifth
reflection. The set of elements $(x,\vphi)$ such that $x^+ \not\in
I$ is designated by $A^{(3)}_r$.

Let us introduce some auxiliary notation. Mark by letter $L$ the
point $x = (\xi, 0)$, and by letter $L^+$, the point $x^+ = (\xi^+,
0)$. Mark by letters $A$ and $B'$, respectively, the centers of
those $i$- and $j$-reflectors that contain the reflection points.
The first of them is a $(\vphi_1,\vphi_2,\del)$-reflector, and the
second, a $(\psi_1,\psi_2,\del)$-reflector; thereby the values
$\vphi_1$,\, $\vphi_2$,\, $\psi_1$,\, $\psi_2$ are determined.
Recall that the angle $\theta$ is such that $\frac{OA}{|OA|} =
e_{\theta}$;\, $\theta^+ := \vphi_\s(\theta)$; and the point $B \in
[P_{j-1},\, P_j]$ is such that $\frac{OB}{|OB|} = e_{\theta^+}$.
Further, define $d\theta^+$ in such a way that $\frac{OB'}{|OB'|} =
e_{\theta^+ + d\theta^+}$;~ define $d\theta$ in such a way that
$\vphi_\s(\theta + d\theta) = \theta^+ + d\theta^+$, and denote by
$A' \in [P_{i-1},\, P_i]$ the point such that $\frac{OA'}{|OA'|} =
e_{\theta + d\theta}$. On the figure below, the case $i < j$,\, $\xi
> 0$ is shown; these relations imply that $d\theta < 0$ and
$d\theta^+ > 0$. The last inequalities are used in the geometric
reasoning given below. The cases where $i > j$ and/or $\xi < 0$ can
be considered analogously.

One easily sees that the values $\vphi_1$ and $\vphi_2$, defined
above, are the angles the rays $AO$ and $AB$, respectively, form
with the normal $-n_i$, and the values $\psi_1$ and $\psi_2$ are the
angles the rays $B'O$ and $B'A'$, respectively, form with $-n_j$.
Remind that throughout this paper, the angles are counted clockwise
from the normal vector to the ray. Therefore, on the figure below
one has $\vphi_1
> 0$,\, $\vphi_2 < 0$,\, $\psi_1 < 0$,\, $\psi_2
> 0$.
 %\newpage

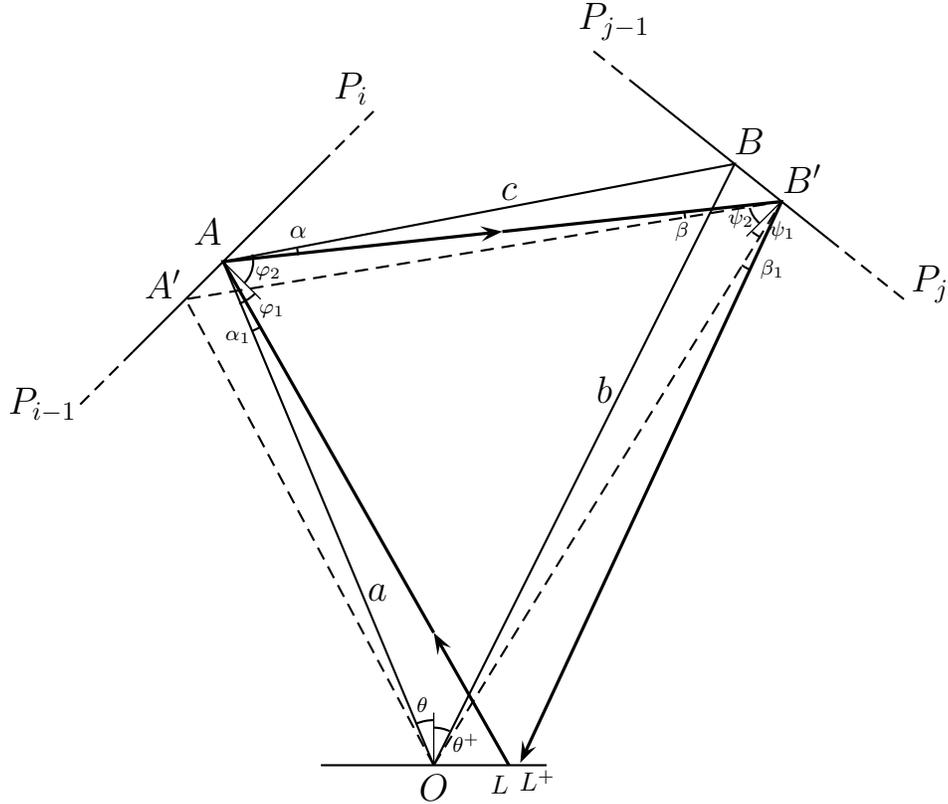
\vspace*{105mm}

\rput(5.7,0){
 \psline(-1.5,0)(1.5,0)
 \psline(0,0)(-2.8,6.7)(4,8)(0,0)
 \psline[linestyle=dashed,linewidth=0.8pt](0,0)(-3.3,6.2)(4.625,7.5)(0,0)
 \psline[linewidth=1.2pt,arrows=->,arrowscale=1.5](1,0)(0,1.763)
 \psline[linewidth=1.2pt,arrows=->,arrowscale=1.5](0,1.763)(-2.8,6.7)(0.9125,7.1)
 \psline[linewidth=1.2pt](0.9125,7.1)(4.625,7.5)
 \psline[linewidth=1.2pt,arrows=->,arrowscale=1.5](4.625,7.5)(1.15,0.03)
      \psline(-4,5.5)(-1.5,8)
    \psline[linestyle=dashed](-4,5.5)(-4.7,4.8)
    \psline[linestyle=dashed](-0.8,8.7)(-1.5,8)
  \rput(-1.1,9){\large $P_i$}
  \rput(-5.2,4.8){\large $P_{i-1}$}
      \psline(5.25,7)(2.75,9)
    \psline[linestyle=dashed](2.125,9.5)(2.75,9)
   \psline[linestyle=dashed](5.25,7)(6.25,6.2)
  \rput(6.6,6.4){\large $P_j$}
  \rput(2.4,9.9){\large $P_{j-1}$}
     \rput(-3,7.1){\large $A$}
     \rput(0,-0.3){\large $O$}
     \rput(4.2,8.3){\large $B$}
     \rput(4.9,7.8){\large $B'$}
     \rput(-3.6,6.4){\large $A'$}
  \rput(1,7.63){\large $c$}
  \rput(2.27,5){\large $b$}
  \rput(-0.75,2.3){\large $a$}
 \psline[linewidth=0.5pt](0,0)(0,0.7)
 \psline[linewidth=0.5pt](-2.8,6.7)(-2.3,6.2)
         \psarc(-2.8,6.7){1.0}{5}{10}
         \psarc(-2.8,6.7){0.4}{-45}{14}
         \psarc(-2.8,6.7){0.6}{-69}{-45}
   \rput(-1.8,7.1){\footnotesize $\al$}
   \rput(-2.2,6.55){\scriptsize $\vphi_2$}
   \rput(-2.15,6.05){\scriptsize $\vphi_1$}
                               \rput(0.52,-0.48){
                               \psline[linewidth=0.5pt](4.05,8)(3.57,7.52)
                               \psarc(4,8){0.4}{-165}{-130}
                               \psarc(4,8){0.55}{-130}{-116}
                               \rput(3.5,7.75){\scriptsize $\psi_2$}
                               \rput(4.05,7.6){\scriptsize $\psi_1$}
                               }
         \psarc(0,0){0.6}{90}{112}
         \psarc(0,0){0.5}{65}{90}
   \rput(-0.15,0.8){\scriptsize $\theta$}
   \rput(0.43,0.3){\scriptsize $\theta^+$}
 \rput(0.88,-0.25){\footnotesize $L$}
 \rput(1.38,-0.2){\footnotesize $L^+$}
         \psarc(-2.8,6.7){1.0}{-67}{-60}
   \rput(-2.6,5.7){\scriptsize $\al_1$}
         \psarc(4.625,7.5){1.3}{-174}{-170}
   \rput(3.3,7.1){\scriptsize $\bt$}
         \psarc(4.625,7.5){1.0}{-122}{-115}
   \rput(4.5,6.6){\scriptsize $\bt_1$}

}

\vspace{17mm}

Designate $a = |OA|$,\, $b = |OB|$,\, $c = |AB|$. Recalling the
definition of the set $A^{(1)}_r$ and taking into account that $|AB|
> |AP_i|$, one concludes that $c \ge 1$. Designate also $dx = |AA'|$,\,
$dx^+ = |BB'|$;\, $\al = \measuredangle BAB'$,\, $\bt =
\measuredangle AB'A'$;\, $\al_1 = \measuredangle OAL$,\, $\bt_1 =
\measuredangle OB'L^+$.

Recall that $|OL| \le 1/2$, the value $a$ is of order $r$, direction
of the ray $LA$ does not coincide with direction of the billiard
trajectory, but rather deviates from it, the difference being
$O(\del/r)$. Thus, $\vphi = \theta + \al_1 + O(\del/r)$. Further,
the difference between the direction of $AB'$ and the true direction
of the billiard particle is $O(\del/c)$. Besides, $\al_1 = O(1/r)$.
The billiard particle intersects the $\kappa(r)$-base of the
$i$-reflector, makes two reflections off its boundary, then
intersects the base again and leaves the reflector, the angles at
the moments of the first and second intersection being
$\tilde{\vphi} = \vphi_1 - \al_1 + O(\del/r)$ and $\tilde{\vphi}^+ =
\vphi_2 + \al + O(\del/c)$, respectively. Applying corollary
\ref{sledstvie c}, one gets $\,\tilde{\vphi}^+ - \vphi_2 =
-\frac{\cos\vphi_1}{\cos\vphi_1}\, (\tilde\vphi - \vphi_1) + O(\del)
+ O((\tilde\vphi - \vphi_1)^2)$, hence
$$
\al + O(\del/c) = -\frac{\cos\vphi_1}{\cos\vphi_2} \left(-\al_1 +
O(\del/r)\right) + O\left(\del\right) +  O\left(1/r^2\right),
$$
therefore
\begin{equation}\label{4a}
\al_1 \cos\vphi_1 = \al \cos\vphi_2 + o(1/r).
\end{equation}
This implies that $\al = O(1/r)$, hence $dx^+ = O(c/r)$ and
$d\theta^+ = O(c/r^2)$. The last equality implies that $d\theta =
O(c/r^2)$, therefore $dx = O(c/r)$ and $\bt = O(1/r)$. After leaving
the $i$-reflector, the particle gets into the $j$-reflector through
its $\kappa(r)$-base, makes the third and fourth reflections off its
boundary, and then gets out of this reflector; the angles at the
moment of getting in and getting out are $\check\vphi = \psi_2 + \bt
+ O(\del/c)$ and $\check{\vphi}^+ = \psi_1 - \bt_1 + O(\del/r)$,
respectively. Applying corollary \ref{sledstvie c} again, one comes
to the formula $\,\check{\vphi}^+ - \psi_1 =
-\frac{\cos\psi_2}{\cos\psi_1}\, (\check\vphi - \psi_2) + O(\del) +
O((\check\vphi - \psi_2)^2)$, hence
$$
-\bt_1 + O(\del/r) = -\frac{\cos\psi_2}{\cos\psi_1} \left(\bt +
O(\del/c)\right) + O\left(\del\right) +  O\left(1/r^2\right),
$$
therefore
\begin{equation}\label{4b}
\bt_1 \cos\psi_1 = \bt \cos\psi_2 + o(1/r).
\end{equation}
The function $\vphi_\s$ preserves the measure $\lam$ and
monotonically decreases on $\Theta_i$, hence
\begin{equation}\label{5}
-\cos\theta d\theta = \cos\theta^+ d\theta^+ (1 + o(1)).
\end{equation}

Finally, simple trigonometric relations for the triangles shown on
the figure imply the following (asymptotic as $r \to +\infty$)
equalities (just for the purpose of completeness, we put below the
formulas (\ref{4a})--(\ref{5})):
$$
\al_1 = \xi\, \frac{\cos\theta}{a}\, (1 + o(1));    \    \   \   \
\al = \al_1\, \frac{\cos\vphi_1}{\cos\vphi_2} + o(1/r); \ \ \ \ dx^+
= \al\, \frac{c}{\cos\psi_2}\, (1 + o(1));
$$
$$
d\theta^+ = dx^+\, \frac{\cos\psi_1}{b}\, (1 + o(1));   \  \  \
-d\theta = d\theta^+\,\frac{\cos\theta^+}{\cos\theta}\,(1 + o(1)); \
\ \ dx = -d\theta\, \frac{a}{\cos\vphi_1}\, (1 + o(1));
$$
$$
\bt = dx\, \frac{\cos\vphi_2}{c}\, (1 + o(1));   \   \   \    \
\bt_1 = \bt\, \frac{\cos\psi_2}{\cos\psi_1} + o(1/r);  \  \ \ \
\xi^+ = \bt_1\, \frac{b}{\cos\theta^+}\, (1 + o(1)).
$$
Making successive substitutions (from the bottom to the top) in
these formulas, one obtains
\begin{equation*}\label{fla}
\xi^+ = \xi(1 + o(1)) + \frac{o(1)}{\cos\theta} +
\frac{o(1)}{\cos\theta^+}\,.
\end{equation*}
This relation can be rewritten in the form
\begin{equation*}
|\xi^+| \le |\xi|\, (1 + \al(r)) + \frac{\al(r)}{\cos\theta} +
\frac{\al(r)}{\cos\theta^+}\,,
\end{equation*}
where $\al(r)$ is a positive function going to zero as
$r\to+\infty$. Taking into account that $\vphi = \theta + O(1/r)$,\
$\vphi^+ = \theta^+ + O(1/r)$, select a constant $c$ such that
$|\vphi - \theta| \le c/r$,\ $|\vphi^+ - \theta^+| \le c/r$. Define
the three sets $A'_r$,\, $A''_r$,\, $A'''_r \subset (I \times
\Theta_i) \setminus (A^{(1)}_r \cup A^{(2)}_r)$ as follows:~ $A'_r$
is the set of pairs $x = (\xi, 0),\, \vphi$ such that $|\xi| \ge
(1/2 - 2\sqrt{\al(r)})/(1 + \al(r))$;~ $A''_r$ is the set of
$(x,\vphi)$ such that $|\vphi| \ge \arccos \sqrt{\al(r)} - c/r$;~
$A'''_r$ is the set of $(x,\vphi)$ such that $|\vphi^+| \ge \arccos
\sqrt{\al(r)} - c/r$. It is not difficult to check that the measure
of each of these sets goes to zero as $r \to +\infty$. Moreover, if
$(x,\vphi)$ does not belong to these sets then $|\xi^+| < 1/2$,
therefore $x^+ \in I$. This implies that $A_r^{(3)} \subset A'_r
\cup A''_r \cup A'''_r$ and therefore $\lim_{r\to+\infty} \mu \big(
A^{(3)}_r \big) = 0$.

Denote $A_r = A^{(1)}_r \cup A^{(2)}_r \cup A^{(3)}_r$; one has
$\lim_{r\to+\infty} \mu(A_r) = 0$. If $(x,\vphi) \in (I \times
\Theta_i) \setminus A_r$ then there are four reflections at points
$\pl\Om$ that do not belong to $I$, and the fifth reflection, at a
point of $I$; therefore, one has $\vphi^+ =
\vphi^+_{\Om(r),I}(x,\vphi)$ and $n_{\Om(r),I}(x,\vphi) = 5$. Thus,
(\ref{convergence n}) is proved.

It follows from the equations $\vphi^+ = \theta^+ + o(1)$,\,
$\theta^+ = \vphi_\s(\theta)$,\, $\theta = \vphi + o(1)$ that
$\vphi^+ = \vphi_\s(\vphi) + o(1)$,\, $r \to +\infty$; therefore,
(\ref{convergence fi}) is also proved.
 \end{proof}

From lemma \ref{l4} it follows

\begin{corollary}
As $r \to +\infty$,\ $\nu_{\Om(r),I}$ weakly converges to $\nu^\s$.
\end{corollary}

\begin{proof}

Indeed, (\ref{conv fi}) implies that the mapping $(\vphi,
\vphi^+_{\Om(r),I}(x,\vphi))$ from $I \times \interval$ to $Q$, as
$r \to +\infty$, converges in measure to the mapping $(\vphi,
\vphi_\s(\vphi))$, hence $f(\vphi, \vphi^+_{\Om(r),I} (x,\vphi))$
converges in measure to $f(\vphi, \vphi_\s(\vphi))$ for any
continuous function $f:\, Q \to \RRR$. This implies that
\begin{equation*}\label{conv f}
\lim_{r \to +\infty}\,\ \int\!\!\!\!\!\!\!\!\!\!\!\!\int\limits_{I
\times \interval} \!\! f(\vphi, \vphi^+_{\Om(r),I} (x,\vphi))\,
d\mu(x,\vphi) =\,\ \int\!\!\!\!\!\!\!\!\!\!\!\!\int\limits_{I \times
\interval} \!\! f(\vphi, \vphi_\s(\vphi))\, d\mu(x,\vphi).
\end{equation*}
Using the relations $\int\!\!\int_{I \times \interval} f(\vphi,
\vphi^+_{\Om,I} (x,\vphi))\, d\mu(x,\vphi) = \int\!\!\!\int_Q
f(\vphi, \vphi^+)\, d\nu_{\Om,I}(\vphi, \vphi^+)$ and
$\int\!\!\int_{I \times \interval} f(\vphi, \vphi_\s(\vphi))\,
d\mu(x,\vphi) = \int\!\!\!\int_Q f(\vphi, \vphi^+)\, d\nu^\s (\vphi,
\vphi^+)$, which arise from the definitions of measures
$\nu_{\Om,I}$,\, $\nu^\s$ and hold true for arbitrary $f$ and $(\Om,
I) \in \SSS$, one gets
\begin{equation*}
\lim_{r \to +\infty} \int_Q f\, d\nu_{\Om(r),I} = \int_Q f\,
d\nu^\s.
\end{equation*}

Corollary is proved.

\end{proof}

According to lemma 2, any measure from $\MMM$ can be approximated by
measures of the kind $\nu^\s$, where $\s^2 = $id,\, $\s(1) \ne m$,
and according to corollary 2, any such measure $\nu^\s$ can be
approximated by measures of the kind $\nu_{\Om(r),I}$ where
$(\Om(r),I) \in \SSS$. This proves theorem 2.

\section*{Appendix}

\subsection*{Proof of lemma \ref{lem vspomogat}}

Let $(\Om, I) \in \SSS$. Taking into account that the measure
$\nu_{\Om,I}$ is invariant with respect to isometries applied
simultaneously to $\Om$ and to $I$, one can assume, without loss of
generality, that $I$ coincides with the segment $[-a,\, a] \times \{
0 \}$,\, $a > 0$, and $\Om \subset \RRR \times [0,\, +\infty)$.
Define the sets $\Om^{[n]}$,\, $I^{[n]}$ by
$$
\Om^{[n]} = \Om \cup \left( [-a,\, a] \times [-1/n,\, 0] \right), \
\ \ \ I^{[n]} = [-a,\, a] \times \{ -1/n \}.
$$

Define the mapping $\TTT_{(n)}$ in the following way. Consider
billiard in the rectangle $[-a,\, a] \times [-1/n,\, 0]$. From a
point $x \in I_n$ let out a particle with initial velocity
$(\sin\vphi, \cos\vphi)$,\, $-\pi/2 < \vphi < \pi/2$. After several
(maybe none) reflections off vertical sides of the rectangle, the
particle reflects off $I$ at a point $x^{(n)}(x,\vphi)$, its
velocity before the reflection being $(\sin\vphi^{(n)}(x,\vphi),\,
\cos\vphi^{(n)}(x,\vphi))$. Thereby, the mapping $\TTT_{(n)}:\,
(x,\vphi) \mapsto (x^{(n)}(x,\vphi),\, \vphi^{(n)}(x,\vphi))$ is
determined; it is defined and takes values on full measure subsets
of $I^{[n]} \times \interval$ and preserve the measure $\mu$. It is
not difficult to see that the mappings $\TTT_{\Om^{[n]},I^{[n]}}$
and $\TTT_{\Om,I}$ are interconnected in the following way:
$$
\TTT_{\Om^{[n]},I^{[n]}} = \TTT_{(n)}^{-1}\, \TTT_{\Om,I}\,
\TTT_{(n)}\,.
$$

Let $B_n = \{ (x,\vphi) \in I^{[n]} \times \interval:\,
\vphi^{(n)}(x,\vphi) = \vphi \}$. It is easy to see that $B_n
\subset \left\{ (x,\vphi):\, x = (\xi, -1/n),\ -a \le \xi \le a,\
-\pi/2 \le \vphi \le -\arctan n(a + \xi) \right.$ or $\left. \arctan
n(a - \xi) \le \vphi \le \pi/2 \right\}$; this implies that
$\mu(B_n) = o(1)$,\, $n \to \infty$, where $d\mu(x,\vphi) =
\cos\vphi\, dx d\vphi$. Put $B_n' =
\TTT_{\Om^{[n]},I^{[n]}}^{-1}(B_n)$.

Let $A \subset Q = \interval \times \interval$ be an arbitrary
measurable set. Denote $D_A = \{ (x,\vphi):\, (\vphi,
\vphi_{\Om,I}^+(x,\vphi)) \in A \}$ and $D_A^{(n)} = \{ (x,\vphi):\,
(\vphi, \vphi_{\Om^{[n]}, I^{[n]}}^+(x,\vphi)) \in A \}$. One has
$$
D_A^{(n)}\, \triangle \left( \TTT_{(n)}^{-1} \left( D_A \right)
\right) \subset B_n \cup B_n',
$$
where $\triangle$ means symmetric difference of sets. Hence
$|\mu(D_A^{(n)}) - \mu\left(D_A\right)| \le \mu\left(B_n \cup
B_n'\right) \le 2\mu\left(B_n\right)$; this implies
$$
\nu_{\Om^{[n]},I^{[n]}}(A) = \nu_{\Om,I}(A) + o(1),
$$
the estimate $o(1)$ being uniform over all $A$. This means that the
sequence of measures $\nu_{\Om^{[n]},I^{[n]}}$ converges in
variation, and hence weakly, to $\nu_{\Om,I}$.

\subsection{Proof of lemma \ref{dense}}

Let us previously prove an auxiliary statement.

\begin{lem}\label{discretization}
Let $A = (a_{ij})_{i,j=1}^m$ be a symmetric matrix,\, $a_{ij}$ being
nonnegative integers. Denote $n_i = \sum_{j=1}^m a_{ij}$. There
exist matrices $B_{ij} = (b_{ij}^{\mu\nu})_{\mu,\nu}$,\, $i,\ j = 1,
\ldots, m$ of order $n_i \times n_j$ such that $B_{ij}^T =
B_{ji}$,\, $\sum_{\mu=1}^{n_i} \sum_{\nu=1}^{n_j} b_{ij}^{\mu\nu} =
a_{ij}$ and the block matrix $B = (B_{ij})$ contains precisely one
''1'' in each row and each column, other elements of $B$ being
zeros.
\end{lem}

\begin{zam}It may happen that for some values $i = i_1,\ i_2,
\ldots$ and for any $j$, holds $a_{ij} = 0$. For these values of
$i$, one has $n_i = 0$ and the matrices $B_{ij}$ and $B_{ji}$ for
any $j$ have the order $0 \times n_j$ and $n_j \times 0$,
respectively, that is, they are empty sets. In this case the matrix
$B$ coincides with the block matrix $B' = (B_{ij})$, where the rows
with indices $i = i_1,\ i_2, \ldots$ and columns with indices $j =
i_1,\ i_2, \ldots$ are crossed out.
\end{zam}

\begin{proof}
Apply induction by $m$. Let the lemma be proved for $m-1$. When
applying it to the matrix $\tilde A = (a_{ij})_{i,j= 2}^m$, one
obtains that there exists a block matrix $\tilde B = (\tilde
B_{ij})_{i,j=2}^m$ satisfying the condition of lemma. Note that the
order of $\tilde B_{ij}$ is $\tilde n_i \times \tilde n_j$, where
$\tilde n_i = \sum_{j=2}^m a_{ij} = n_i - a_{i1}$. Define the
matrices $B_{ij}$ as follows.

(a) Let $B_{11} =$ diag\,$\{ \underbrace{1, \ldots, 1}_{a_{11}},\ 0,
\ldots, 0 \}$.

(b) Let $b_{12}^{a_{11}+1,1} = \ldots =
b_{12}^{a_{11}+a_{12},a_{12}} = 1$;\ $b_{13}^{a_{11}+a_{12}+1,1} =
\ldots = b_{13}^{a_{11}+a_{12}+a_{13},a_{13}} = 1$;\ \ldots; \
$b_{1m}^{a_{11}+\ldots+a_{1,m-1}+1,1} = \ldots =
b_{1m}^{a_{11}+\ldots+a_{1m},a_{1m}} = 1$;\ other elements of the
matrices $B_{1j}$,\, $j = 2, \ldots, m$ are zeros. Thus, in each
matrix $B_{1j}$, on the diagonal, whose first element belongs to the
first column and the $(a_{11} + a_{12} + \ldots + a_{1,j-1} + 1)$th
row, the first $a_{1j}$ elements equal 1, and the remaining elements
on this diagonal and all the elements outside the diagonal equal
aero. Thus, the matrices $B_{1j}$,\, $j = 2, \ldots, m$ are defined.
The matrices $B_{i1}$,\, $i = 2, \ldots, m$ are defined by the
condition $B_{i1} = B_{1i}^T$.

(c) For $i \ge 2$,\, $j \ge 2$ define the matrix $B_{ij}$ as
follows. For $\mu \le a_{1i}$ or $\nu \le a_{1j}$ let
$b_{ij}^{\mu\nu} = 0$, and for $\mu \ge a_{1i} + 1$,\, $\nu \ge
a_{1j} + 1$ let $b_{ij}^{\mu\nu} = \tilde
b_{ij}^{\mu-a_{1i},\nu-a_{1j}}$. Thus, in the obtained matrix
$B_{ij}$, the right lower corner coincides with the matrix $\tilde
B_{ij}$, and all the remaining elements are equal to zero. The
number of rows of this matrix equals $a_{1i} + \tilde n_i = n_i$,
and the number of columns, $a_{1j} + \tilde n_j = n_j$. It is also
not difficult to verify that $\sum_{\mu\nu} b_{ij}^{\mu\nu} =
a_{ij}$ and that each row and each column of the obtained block
matrix $B = (B_{ij})_{i,j=1}^m$ contains precisely one number ''1''.
\end{proof}

Now, let us proceed to the proof of lemma \ref{dense}.

Consider an arbitrary measure $\nu \in \MMM$. Let $n \in \NNN$;
consider the partition of the square $Q$ to $n^2$ smaller squares
$Q_{ij}^n = \Theta_i^n \times \Theta_j^n$;\, $Q = \cup_{i,j=1}^n
Q_{ij}^n$. Let us show that for all $n \in \NNN$,\, $1 \le i,\, j
\le n$ one can choose rational non-negative numbers $c_{ij}^n$ in
such a way that $c_{ij}^n = c_{ji}^n$, for any $i \ \, \sum_j
c_{ij}^n = 2/n$ and $n^2 \cdot \max_{i,j} |\nu\left( Q_{ij}^n
\right) - c_{ij}^n|$ tends to 0 as $n \to \infty$.

Indeed, for $i > j$ choose values $c_{ij}^n$ in such a way that
$|c_{ij}^n - \nu(Q_{ij}^n)| \le n^{-4}$;~ for $i < j$ put $c_{ij}^n
= c_{ji}^n$;~ finally, for $i = j$ put $c_{ii}^n = 2/n - \sum_{j\ne
i} c_{ij}^n$. One has $\sum_{j=1}^n \nu(Q_{ij}^n) =
\nu\left(\Theta_i^n \times \interval\right) = \lam(\Theta_i^n) =
2/n$, hence $\nu(Q_{ii}^n) = 2/n - \sum_{j\ne i} \nu(Q_{ij}^n)$.
This implies that $|c_{ii}^n - \nu(Q_{ii}^n)| = \big|\sum_{j\ne i}
\left( c_{ij}^n - \nu(Q_{ij}^n) \right)\big| \le (n - 1) \cdot
n^{-4} < n^{-3}$. Therefore $|\max_{i,j} \left( \nu(Q_{ij}^n \right)
- c_{ij}^n)| < n^{-3}$.

Any sequence of measures $\nu_n$, satisfying the conditions
$\nu_n(Q_{ij}^n) = c_{ij}^n$,\,\, $1 \le i,\, j \le n$, weakly
converges to $\nu$. Indeed, for any continuous function $f:\, Q \to
\RRR$ holds
$$
\int_Q f\, d\nu_n - \int_Q f\, d\nu = \sum_{i,j=1}^n \left(
\int_{Q_{ij}^n} f\, d\nu_n - \int_{Q_{ij}^n} f\, d\nu \right) =
$$
$$
= \sum_{i,j=1}^n \left( f(x_{ij}^n)\, c_{ij}^n - f(\tilde
x_{ij}^n)\, \nu\left( Q_{ij}^n \right) \right) =
$$
\begin{equation}\label{abc}
= \sum_{i,j=1}^n f(x_{ij}^n) \left( c_{ij}^n - \nu\left( Q_{ij}^n
\right) \right) + \sum_{i,j=1}^n \left( f(x_{ij}^n) - f(\tilde
x_{ij}^n) \right) \nu\left( Q_{ij}^n \right),
\end{equation}
where $x_{ij}^n,\ \tilde x_{ij}^n \in Q_{ij}^n$. Taking into account
that $|f(x_{ij}^n) - f(\tilde x_{ij}^n)| \le \max_{i,j}
\sup_{x',x''\in Q_{ij}^n} |f(x') - f(x'')| \to 0$ as $n \to \infty$
and that $n^2 \cdot \max_{i,j} |\nu \left( Q_{ij}^n \right) -
c_{ij}^n| \to 0$ as $n \to \infty$, one concludes that the
expression in (\ref{abc}) tends to zero as $n \to \infty$.

Thus, to finish the proof of lemma, it suffices to show that for any
$n \in \NNN$ there exist a positive integer $m$  and a permutation
$\s_n$ of $\{ 1, \ldots, m \}$ such that $\s_n^2 =$ id,\,
$\nu^{\s_n}(Q_{ij}^n) = c_{ij}^n$ and $\s(1) \ne m$.

Fix $n$ and find a positive integer $N$ such that all values $a_{ij}
:= N \cdot c_{ij}^n$ are integer. Obviously, the matrix $A =
(a_{ij})_{i,j=1}^n$ is symmetric and for any $i\,$\ $\sum_{j=1}^n
a_{ij} = 2N/n \in \NNN$. According to lemma \ref{discretization},
there exist square matrices $B_{ij}$ of order $2N/n$ such that
$B_{ij}^T = B_{ji}$, the sum of elements of each matrix $B_{ij}$
equals $a_{ij}$, and composed of them block matrix $B = (B_{ij})$
has precisely one unit in each row and each column, and all the
remaining elements are zeros.

Consider block matrices $D_{ij}$ of double size, obtained from
$B_{ij}$ by substitution of elements with $2 \times 2$ matrices:
zeros are substituted with $\left(
\begin{array}{cc} 0 & 0\\ 0 & 0
\end{array} \right)$, and units, with $\left( \begin{array}{cc} 1 & 0\\ 0 &
1 \end{array} \right)$. The matrices $D_{ij}$ have order $4N/n$,
holds $D_{ij}^T = D_{ji}$, the sum of elements in each matrix
$D_{ij}$ equals $2a_{ij}$, and the block matrix $D = (D_{ij})$
composed of $D_{ij}$ has precisely one unit in each column and each
row, and all remaining elements are zeros.

The order of $D$ is $4N$. Let $m = 4N$ and denote by $d_{ij}$,\,
$i,\, j = 1, \ldots, m$, the elements of $D$. The peculiarity of
$D$, as compared to $B$, is that the equality $d_{1m} = 0$ is
guaranteed. Define the mapping $\s = \s_n:\, \{ 1, \ldots, m \} \to
\{ 1, \ldots, m \}$ in such a way that for any $i, \ d_{i\s(i)} =
1$. The so defined mapping $\s_n$ is a permutation; since the matrix
$D$ is not symmetric, one concludes that $\s_n^2 =$ id, and since
$d_{1m} \ne 1$, one concludes that $\s_n(1) \ne m$. Besides, the
following holds true: $\nu^{\s_n}(Q_{ij}^n) = N^{-1} \sum_{\mu,\nu}
b_{ij}^{\mu\nu} = c_{ij}^n$. Lemma is proved.

\section*{Acknowledgements}

This work was supported by {\it Centre for Research on Optimization
and Control} (CEOC) from the ''{\it Funda\c{c}\~{a}o para a
Ci\^{e}ncia e a Tecnologia}'' (FCT), cofinanced by the European
Community Fund FEDER/POCTI.


\begin{thebibliography}{99}

\bibitem{N}
I. Newton,\, {\it Philosophiae naturalis principia mathematica}\,
1686.

\bibitem{BK} G. Buttazzo and B. Kawohl.
\textit{On Newton's problem of minimal resistance}. Math. Intell.
{\bf 15}, 7--12 (1993).

\bibitem{BrFK} F. Brock, V. Ferone, and B. Kawohl.
\textit{A symmetry problem in the calculus of variations}. Calc.
Var. {\bf 4}, 593--599 (1996).

\bibitem{BFK} G. Buttazzo, V. Ferone, and B. Kawohl.
\textit{Minimum problems over sets of concave functions and related
questions}. Math. Nachr. {\bf 173}, 71--89 (1995).

\bibitem{BG}
G. Buttazzo and P. Guasoni,\, {\it Shape optimization problems over
classes of convex domains},\, J. Convex Anal. {\bf 4}, No.2, 343-351
(1997).

\bibitem{LP1} T. Lachand-Robert and M.~A. Peletier.
\textit{Newton's problem of the body of minimal resistance in the
class of convex developable functions}. Math. Nachr. {\bf 226},
153--176 (2001).

\bibitem{LP2} T. Lachand-Robert, M.~A. Peletier.
\textit{An example of non-convex minimization and  an application to
Newton's problem of the body of least resistance}. Ann. Inst. H.
Poincar\'e, Anal. Non Lin. {\bf 18}, 179--198 (2001).

\bibitem{CL1} M. Comte and T. Lachand-Robert.
\textit{Newton's problem of the body of minimal resistance under a
single-impact assumption}. Calc. Var. Partial Differ. Equ. {\bf 12},
173--211 (2001).

\bibitem{CL2} M. Comte and T. Lachand-Robert.
\textit{Existence of minimizers for Newton's problem of the body of
minimal resistance under a single-impact assumption}. J. Anal. Math.
{\bf 83}, 313--335 (2001).

\bibitem{LO} T. Lachand-Robert and E. Oudet.
\textit{Minimizing within convex bodies using a convex hull method}.
SIAM J. Optim. {\bf 16}, 368--379 (2006).

\bibitem{P1}
A.\,Yu. Plakhov. {\it Newton's problem of a body of minimal
aerodynamic resistance}, Dokl. Akad. Nauk {\bf 390}, N$^\text{o}$3,
314--317 (2003).

\bibitem{P2}
A.\,Yu. Plakhov. {\it Newton's problem of the body of minimal
resistance with a bounded number of collisions}, Russ. Math. Surv.
{\bf 58} N$^\text{o}$1, 191-192 (2003).

\bibitem{sb-math-averaged04}
A.\,Yu. Plakhov. {\it Newton's problem of the body of minimum mean
resistance}, Sbornik: Mathematics {\bf 195}, N$^\text{o}$7-8,
1017-1037 (2004).

\bibitem{temperature}
A.\,Yu. Plakhov and Delfim F.~M.~Torres. {\it Newton's aerodynamic
problem in media of chaotically moving particles}, Sbornik:
Mathematics {\bf 196}, N$^\text{o}$5-6, 885-933 (2005).

\bibitem{P-unbounded regions}
A.\,Yu. Plakhov. {\it Billiards in unbounded domains reversing the
direction of motion of a particle}, Russ. Math. Surv. {\bf 61}
N$^\text{o}$1, 179-180 (2006).

\bibitem{T}
S. Tabachnikov,\, {Billiards},\,  {\it Paris: Soci\'et\'e\
Math\'ematique de France} (1995).


\end{thebibliography}
\end{document}